\documentclass[final,secthm,amsthm]{elsart}
\usepackage{amssymb,amstext,amsmath,xy}
\xyoption{all}
\def\MSC{\par\leavevmode\hbox {\it 2000 MSC:\ }}%

\newcommand{\bx}{\hfill$\square$}
  \newcommand{\A}{{\mathcal{A}}}
  \newcommand{\B}{{\mathcal{B}}}
  \newcommand{\C}{{\mathcal{C}}}

  \newcommand{\F}{{\mathcal{F}}}
  \newcommand{\G}{{\mathcal{G}}}
\renewcommand{\H}{{\mathcal{H}}}

  \newcommand{\K}{{\mathcal{K}}}

\newcommand{\M}{{\mathcal{M}}}
  
\renewcommand{\O}{{\mathcal{O}}}

  \newcommand{\U}{{\mathcal{U}}}
  \newcommand{\V}{{\mathcal{V}}}
  \newcommand{\W}{{\mathcal{W}}}

\newcommand{\bC}{{\mathbb{C}}}

\newcommand{\bF}{{\mathbb{F}}}
\newcommand{\bN}{{\mathbb{N}}}

\newcommand{\bT}{{\mathbb{T}}}
\newcommand{\bZ}{{\mathbb{Z}}}
\newcommand{\fg}{{\mathfrak{g}}}

\renewcommand{\phi}{\varphi}
\newcommand{\qand}{\quad\text{and}\quad}
\newcommand{\qif}{\quad\text{if}\quad}
\newcommand{\qfor}{\quad\text{for}\quad}
\newcommand{\qforal}{\quad\text{for all}\quad}
\newcommand{\AnD}{\text{ and }}
\newcommand{\FOR}{\text{ for }}
\newcommand{\FORAL}{\text{ for all }}
\newcommand{\OR}{\text{ or }}
\newcommand{\Ad}{\operatorname{Ad}}
\newcommand{\diag}{\operatorname{diag}}
\newcommand{\spn}{\operatorname{span}}
\DeclareMathOperator*{\wotlim}{\textsc{wot}--lim}

\newcommand{\bsl}{\setminus}
\newcommand{\ca}{\mathrm{C}^*}
\newcommand{\Fn}{\mathbb{F}_n^+}
\newcommand{\Fm}{\mathbb{F}_m^+}
\newcommand{\Fth}{\mathbb{F}_\theta^+}
\newcommand{\ip}[1]{\langle #1 \rangle}
\newcommand{\mt}{\varnothing}
\newcommand{\ol}{\overline}
\newcommand{\ltwo}{\ell^2}

\begin{document}
\begin{frontmatter}
\title{Atomic Representations\\ of Rank 2 Graph Algebras}

\author[UW]{Kenneth R. Davidson},
\ead{krdavids@uwaterloo.ca}
\thanks{First author partially supported by an NSERC grant.}
\author[Lan]{Stephen C. Power},
\ead{s.power@lancaster.ac.uk}
\thanks{Second author partially supported by EPSRC grant EP/E002625/1.}
\author[UW]{Dilian Yang}
\ead{\\dyang@uwaterloo.ca}
\thanks{Third author partially supported by the Fields Institute.}

\address[UW]{Pure Math. Dept.,U. Waterloo,
Waterloo, ON  N2L--3G1,CANADA}
\address[Lan]{Dept. Maths. Stats., Lancaster University,
Lancaster LA1 4YF, U.K. }

\begin{abstract}
We provide a detailed analysis of atomic $*$-repre\-sent\-ations
of rank 2 graphs on a single vertex.
They are completely classified up to unitary
equivalence, and decomposed into a direct sum or direct
integral of irreducible atomic representations.
The building blocks are described as the minimal $*$-dilations
of defect free representations modelled on finite groups of rank 2.
\end{abstract}


\end{frontmatter}

\section{Introduction}\label{S:intro}

Kumjian and Pask \cite{KumPask} introduced higher rank graphs
and their associated C*-algebras as a generalization of graph C*-algebras
that are related to the generalized Cuntz--Kreiger algebras of
Robertson and Steger \cite{RobSims}.
The C*-algebras of higher rank graphs have been studied in a variety of papers
\cite{RaeSimYee1,RaeSimYee2,FarMuhYee,PaskRRS,Sims}.
See also Raeburn's survey \cite{Raeburn}.

In \cite{KP1}, Kribs and Power examined the nonself-adjoint
operator algebras which are associated with these higher rank graphs.
More recently, Power \cite{P1} presented a detailed analysis
of single vertex rank 2 case.
As this case already contains many new and intriguing
algebras, we were motivated to study them more closely.

The rank 2 graphs on one vertex form an intriguing family
of semigroups with a rich representation theory.
There are already many interesting and non-trivial issues.
Our purpose in this paper is to completely classify the atomic
$*$-rep\-re\-sen\-ta\-tions of these semigroups.
These representations combine analysis with some interesting
combinatorial considerations.
They provide a rich class of representations of the associated
C*-algebra which have proven effective in understanding the
underlying structure.

The algebras are given concretely in terms of a finite set of generators
and relations of a special type.
Given a permutation $\theta$ of $m \times n$,  form a unital
semigroup $\Fth$ with generators $e_1,\dots,e_m$ and $f_1,\dots,f_n$
which is free in the $e_i$'s and free in the $f_j$'s, and has the
commutation relations $e_i f_j = f_{j'} e_{i'}$ where $\theta(i,j)
= (i',j')$ for $1 \le i \le m$ and $1 \le j \le n$. This is a
cancellative semigroup with unique factorization \cite{KumPask}.

A $*$-rep\-re\-sen\-ta\-tion of the semigroup $\Fth$ is a representation
$\pi$  of $\Fth$ as isometries with the property that
\[
 \sum_{i=1}^m \pi(e_i)\pi(e_i)^* = I
 = \sum_{j=1}^n \pi(f_j) \pi(f_j)^* .
\]
An \textit{atomic} $*$-rep\-re\-sen\-ta\-tion acts on a Hilbert space
with an orthonormal basis which is permuted, up to
unimodular scalars, by each of the generators.
The C*-algebra $\ca(\Fth)$ is the universal C*-algebra generated
by these  $*$-rep\-re\-sen\-ta\-tions.

The motivation for studying these representations comes
from the case of the free semigroup $\Fm$ generated by
$e_1,\dots,e_m$ with no relations.
The C*-algebra $\ca(\Fm)$ is just the Cuntz algebra $\O_m$ \cite{Cun}
(see also \cite{KRDCalg}).
Davidson and Pitts \cite{DP1} classified the atomic
$*$-rep\-re\-sen\-ta\-tions of $\Fm$ and showed that the irreducibles fall
into two types, known as ring representations and infinite tail
representations.
They provide an interesting class of C*-algebra
representations of $\O_m$ which are amenable to analysis.
Furthermore the atomic $*$-rep\-re\-sen\-ta\-tions  feature significantly
in the dilation theory of row contractions \cite{DKS}.

The 2-graph situation turns out to be considerably more
complicated than the case of the free semigroup.  Whereas the
irreducible atomic $*$-rep\-re\-sen\-ta\-tions of $\bF_n^+$ are of two
types, the irreducible atomic $*$-rep\-re\-sen\-ta\-tions of $\Fth$ fall
into six types. Nevertheless, we are able to put them all into a
common framework modelled on abelian groups of rank 2.

\section{Background}

\medbreak\textbf{Rank 2 graphs.}
The semigroup $\Fth$ is generated by $e_1,\dots,e_m$ and
$f_1,\dots,f_n$. The identity is denoted as $\mt$. There are no
relations among the $e$'s, so they generate a copy of the free
semigroup on $m$ letters, $\Fm$; and there are no relations on the
$f$'s, so they generate a copy of $\Fn$. There are
\textit{commutation relations} between the $e$'s and $f$'s given
by a permutation $\theta$ in $S_{m\times n}$ of $m \times n$:
\[ e_i f_j = f_{j'} e_{i'} \quad\text{where } \theta(i,j) = (i',j') .\]

A word $w\in\Fth$ has a fixed number of $e$'s and $f$'s regardless
of the factorization; and the \textit{degree} of $w$ is $d(w) := (k,l)$ if
there are $k$ $e$'s and $l$ $f$'s.  The degree map is a homomorphism
of $\Fth$ into $\bN_0^2$. The \textit{length} of $w$ is
$|w| = k+l$. The commutation relations allow any word $w\in\Fth$
to be written with all $e$'s first, or with all $f$'s first, say
$w = e_uf_v = f_{v'}e_{u'}$. Indeed, one can factor $w$ with any
prescribed pattern of $e$'s and $f$'s as long as the degree is
$(k,l)$. It is straightforward to see that the factorization is
uniquely determined by the pattern and that $\Fth$ has the
unique factorization property and cancellation.
See also \cite{KumPask,KP1,P1}.

\begin{exmp}\label{E:favourite}
With $n=m=2$ we note that the relations
\begin{alignat*}{2}
e_1f_1 &= f_2e_1 &\qquad\qquad  e_1f_2 &= f_1e_2 \\
e_2f_1 &= f_1e_1 &\qquad\qquad  e_2f_2 &= f_2e_2
\end{alignat*}
arise from the $3$-cycle permutation
$\theta = \big( (1,1) , (1,2) , (2,1) \big)$
in $S_{2 \times 2}$.
We will refer to $\Fth$ as the forward 3-cycle semigroup.
The reverse $3$-cycle semigroup is the one arising from
the $3$-cycle $\big( (1,1) , (2,1) , (1,2) \big)$.

It was shown by Power in  \cite{P1} that the 24 permutations
of $S_{2 \times 2}$  give rise to $9$ isomorphism
classes of semigroups $\Fth$, where we allow
isomorphisms to exchange the $e_i$'s for $f_j$'s.
In particular, the forward and reverse $3$-cycles give
non-isomorphic semigroups.
\end{exmp}

\begin{exmp}\label{E:flip}
With $n=m=2$  the relations
\begin{alignat*}{2}
e_1f_1 &= f_1e_1 &\qquad\qquad  e_1f_2 &= f_1e_2 \\
e_2f_1 &= f_2e_1 &\qquad\qquad  e_2f_2 &= f_2e_2
\end{alignat*}
are those arising from the $2$-cycle permutation
$\big( (1,2) , (2,1) \big)$ and we refer to $\Fth$
in this case as the flip semigroup because of the
commutation rule: $e_if_j = f_i e_j$ for $1 \le i,j \le 2$.
This is an example which has periodicity,
a concept which will be explained in more detail later.
\end{exmp}

\bigbreak\textbf{Representations.}
We now define two families of representations which will be considered here:
the $*$-rep\-re\-sen\-ta\-tions and defect free (partially isometric) representations.

\begin{defn}
A \textit{partially isometric representation} of $\Fth$ is a
semigroup homomorphism $\sigma: \Fth\to \B(\H)$ whose range
consists of partial isometries on a Hilbert space $\H$.
The representation $\sigma$ is  \textit{isometric} if the range
consists of isometries.
A representation is \textit{atomic} if it is partially isometric and
there is an orthonormal basis which is permuted,
up to scalars, by each partial isometry.
That is, $\sigma $ is atomic if there is a basis  $\{\xi_k : k\ge1\}$
so that for each $w \in \Fth$, $\sigma(w) \xi_k = \alpha \xi_l$ for
some $l$ and some $\alpha \in \bT \cup \{0\}$.

We say that $\sigma$ is \textit{defect free} if
\[
 \sum_{i=1}^m \sigma(e_i) \sigma(e_i)^* = I
 = \sum_{j=1}^n \sigma(f_j) \sigma(f_j)^* .
\]
An isometric defect free representation is called a
\textit{$*$-rep\-re\-sen\-ta\-tion} of $\Fth$.
\end{defn}

For a (partially) isometric representation, the defect free condition
is equivalent to saying that the ranges of the $\sigma(e_i)$'s are
pairwise orthogonal and sum to the identity, and that the same holds
for the ranges of the $\sigma(f_j)$'s.
Equivalently, this says that $[\sigma(e_1)\ \cdots\ \sigma(e_m)]$ is a
(partial) isometry from the direct sum of $m$ copies of $\H$ onto $\H$,
and the likewise for $[\sigma(f_1)\ \cdots\ \sigma(f_n)]$.
A representation which satisfies the property  that these row operators
are isometries \textit{into} $\H$ is called \textit{row isometric}.
The left regular representation $\lambda$ of $\Fth$ is an example
of a representation which is row isometric,  but is not defect free.

\begin{defn}
The C*-algebra $\ca(\Fth)$ is the universal C*-algebra for $*$-rep\-re\-sen\-ta\-tions.
This is the C*-algebra generated by isometries $E_1,\dots,E_m$
and $F_1,\dots,F_n$ which are defect free:
$ \sum_{i=1}^m E_i E_i^* = I
 = \sum_{j=1}^n F_j F_j^*$, and satisfies
the commutation relations of $\Fth$, with the universal property that every
$*$-rep\-re\-sen\-ta\-tion $\sigma$ of $\Fth$ extends to a
$*$-homomorphism of $\ca(\Fth)$ onto $\ca(\sigma(\Fth))$.
\end{defn}

In the case of the free semigroup $\bF_n^+$,
Davidson and Pitts \cite{DP1} classified the atomic
row isometric representations.
They showed that in the irreducible case there are three
possibilities, namely (i) the left regular representation,
(ii) a ring representation, determined by a  primitive word
$u$ in  $\bF_n^+$ and a unimodular scalar $\lambda$,
and (iii) a tail representation, constructed from an aperiodic
infinite word in the generators of $\bF_n^+$.
The left regular representation is the only one which is
not defect free.
The ring and tail representations provide the irreducible atomic
$*$-rep\-re\-sen\-ta\-tions.
The universal C*-algebra of $\Fn$ is the Cuntz algebra $\O_n$.

A ring representation is determined by a set of $k$ basis
vectors  which are cyclically permuted, modulo $\lambda$,
according to the $k$ letters of $u$.
The primitivity of $u$ means that $u$ is not a proper
power of a smaller word. On the other hand, a
tail representation $\sigma$
is determined by an infinite word $z = z_0z_{-1}z_{-2}\dots $ in the
generators of $\bF_n^+$.
There is a subset of basis elements $\xi_0, \xi_{-1}, \xi_{-2},  \dots$
for which $\sigma(z_k)\xi_{k-1} = \xi_k$ for $k \le 0$.
In both cases the subspace $\M$ spanned by the basis vector
subset is coinvariant for $\sigma$ (i.e. $\sigma(w)^*\M \subset \M$ for all $w\in\bF_n^+$)
and cyclic for $\sigma$ (i.e. $\bigvee_{w\in\bF_n^+} \sigma(w)\M = \H$).
On the complementary invariant subspace $\M^\perp$, $\sigma$
decomposes as a direct sum of copies of the left regular
representation.

We shall meet these representations, as restrictions, in the
classification of irreducible atomic representations of $\Fth$.

\medbreak\textbf{Dilations.}
If $\sigma$ is a representation of $\Fth$ on a Hilbert space $\H$,
say that a representation $\pi$ of $\Fth$ on a Hilbert space
$\K \supset \H$ is a \textit{dilation} of $\sigma$ if
\[ \sigma(w) = P_\H \pi(w)|_\H \qforal w \in \Fth .\]
The dilation $\pi$ of $\sigma$ is a \textit{minimal isometric dilation}
if $\pi$ is isometric and the smallest invariant subspace containing
$\H$ is all of $\K$, which means that $\K = \bigvee_{w \in \Fth} \pi(w)\H$.
This minimal isometric dilation is called \textit{unique} if for any two minimal
isometric dilations $\pi_i$ on $\K_i$, there is a unitary operator $U$
from $\K_1$ to $\K_2$ such that $U|_\H$ is the identity map and
$\pi_2 = \Ad U\, \pi_1$.

It is straightforward to see that an isometric dilation of a defect free
representation is still defect free, and hence is a $*$-rep\-re\-sen\-ta\-tion.
An important result from our paper \cite{DPYdiln} is that every
defect free representation has a unique minimal $*$-dilation.
Both existence and uniqueness are of critical importance.
Moreover, if the original representation is atomic, then so is
the $*$-dilation.

\begin{thm}[\cite{DPYdiln}, Theorems~5.1, 5.5] \label{defectfreediln}
Let $\sigma$ be a defect free (atomic) representation.
Then $\sigma$ has a unique minimal
dilation to a (atomic) $*$-rep\-re\-sen\-ta\-tion.
\end{thm}

Conversely, if $\pi$ is a $*$-rep\-re\-sen\-ta\-tion, one can obtain a
defect free representation $\sigma$ on a subspace $\H$ if $\H$ is co-invariant,
i.e. $\H^\perp$ is invariant for $\pi(\Fth)$, by setting $\sigma(w) = P_\H \pi(w)|_\H$.
By the Dilation Theorem~\ref{defectfreediln},
$\sigma$ has a unique minimal $*$-dilation.
But such a $*$-dilation is evidently the restriction of $\pi$ to the
reducing subspace $\K' = \bigvee_{w\in\Fth} \pi(w)\H$.
If the subspace $\H$ is \textit{cyclic}, i.e. $\K'=\K$, then this
minimal $*$-dilation is $\pi$ itself.
In this case, $\pi$ is uniquely determined by $\sigma$.

The significance of this for us is that every atomic $*$-rep\-re\-sen\-ta\-tion
has a particularly nice coinvariant subspace on which the compression
$\sigma$ has a very tractable form.  It is by determining these smaller
defect free representations that we are able to classify the atomic
$*$-rep\-re\-sen\-ta\-tions.

\section{Examples of Atomic Representations}

We begin with a few examples of atomic representations.

\begin{exmp}\label{E:cycle reps}
Given a permutation $\theta$ of $m \times n$,
select a cycle of $\theta$, say
\[ \big( (i_1,j_1), (i_2,j_2), \dots, (i_k,j_k) \big) .\]
Form a Hilbert space of dimension $k$ with
basis $\xi_s$ for $1 \le s \le k$.
Define
\[
 \sigma(e_i) = \sum_{i_s=i} \xi_s \xi_{s-1}^* \qand
 \sigma(f_j) = \sum_{j_s=j} \xi_{s-1} \xi_s^* .
\]
That is, $\sigma(f_{j_s})$ maps $\xi_s$ to $\xi_{s-1}$
and $\sigma(e_{i_s})$ maps $\xi_{s-1}$ back to $\xi_s$.
Likewise $\sigma(e_{i_{s+1}})$ maps $\xi_s$ to $\xi_{s+1}$
and $\sigma(f_{j_{s+1}})$ maps $\xi_{s+1}$ back to $\xi_s$.
This corresponds to the commutation relation
$e_{i_s}f_{j_s} = f_{j_{s+1}}e_{i_{s+1}}$ indicated by
the cycle of $\theta$.
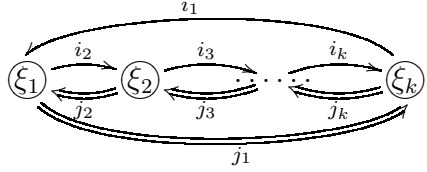
\[
 \xymatrix{
 {\xy *++={\xi_1} *\frm{o} \endxy}
 \ar@/^/[r]^{i_2} \ar@{=>}@(d,dr) [rrr]_{j_1}  
 &{\xy *++={\xi_2} *\frm{o} \endxy}
 \ar@/^/[r]^{i_3} \ar@{=>}@/^/[l]^{j_2}
 &\dots\dots
 \ar@/^/[r]^{i_k} \ar@{=>}@/^/[l]^{j_3}
 &{\xy *++={\xi_k} *\frm{o} \endxy}
 \ar@(u,ul) [lll]_{i_1} \ar@{=>}@/^/[l]^{j_k}\\&&&&
}
\]

So it is not difficult to verify that this defines a defect free
atomic representation of $\Fth$.
By the Dilation Theorem~\ref{defectfreediln}, this can be dilated to a unique
$*$-rep\-re\-sen\-ta\-tion of $\Fth$.

One can further adjust this example by introducing
scalars.  For example, if $\alpha,\beta \in \bT$,
define $\sigma_{\alpha,\beta}(e_i) = \alpha \sigma(e_i)$
and $\sigma_{\alpha,\beta}(f_j) = \beta \sigma(f_j)$.
Two such representations will be shown to be
unitarily equivalent if and only if $\alpha_1^k = \alpha_2^k$
and $\beta_1^k=\beta_2^k$.
\end{exmp}

\begin{exmp}\label{3a_repn} \textbf{Inductive representations.}
An important family of $*$-repre\-sent\-ations were introduced
in \cite{DPYdiln}.
Recall that the left regular representation $\lambda$ acts on $\ltwo(\Fth)$,
which has orthonormal basis $\{\xi_v : v \in \Fth\}$, by left multiplication:
$\lambda(w) \xi_v = \xi_{wv}$ for all $w,v \in \Fth$.

Start with an arbitrary \textit{infinite word} or \textit{tail}
$\tau = e_{i_0}f_{j_0}e_{i_1}f_{j_1} \dots$.
Let $\F_s = \F := \Fth$, for $s=0,1,2,\dots $, viewed as a discrete set
on which the generators of $\Fth$ act as injective maps
by right multiplication, namely,
\[ \rho(w)g = gw \qforal g \in \F. \]
Consider $\rho_s = \rho(e_{i_s}f_{j_s})$ as a map from
$\F_s$ into $\F_{s+1}$.
Define $\F_\tau$ to be the injective limit set
\[
 \F_\tau = \lim_{\rightarrow} (\F_s, \rho_s ) ;
\]
and let $\iota_s$ denote the injections of $\F_s$ into $\F_\tau$.
Thus $\F_\tau$ may be viewed as the union of $\F_0, \F_1, \dots $
with respect to these inclusions.

The left regular action $\lambda$ of $\Fth$ on itself induces
corresponding maps on $\F_s$ by  $\lambda_s(w) g = wg$.
Observe that $\rho_s \lambda_s = \lambda_{s+1} \rho_s$.
The injective limit of these actions is an action $\lambda_\tau$
of $\Fth$ on $\F_\tau$.
Let $\lambda_\tau$ also denote the corresponding representation of $\Fth$
on $\ltwo(\F_\tau)$.
The standard basis of $\ltwo(\F_\tau)$ is  $\{ \xi_g : g \in \F_\tau\}$.
A moment's reflection shows that this provides a defect free,
isometric (atomic) representation of $\Fth$; i.e.\ it is a $*$-rep\-re\-sen\-ta\-tion.

We now describe a coinvariant cyclic subspace
that contains all of the essential information about
this representation.
Let $\H = \ol{ \lambda_\tau(\Fth)^* \xi_{\iota_0(\mt)}}$.
This is coinvariant by construction.
As it contains $\xi_{\iota_s(\mt)}$ for all $s \ge1$,
it is easily seen to be cyclic.
Let $\sigma_\tau$ be the compression of $\lambda_\tau$ to $\H$.
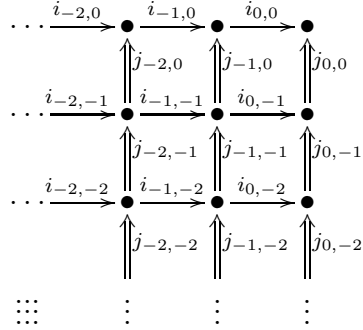
\[
 \xymatrix{
 \dots \ar[r]^{i_{-2,0}} &
 \bullet \ar[r]^{i_{-1,0}} &
 \bullet \ar[r]^{i_{0,0}} &
 \bullet \\
 \dots \ar[r]^{i_{-2,-1}} &
 \bullet \ar[r]^{i_{-1,-1}} \ar@{=>}[u]_(.6){j_{-2,0}} &
 \bullet \ar[r]^{i_{0,-1}} \ar@{=>}[u]_(.6){j_{-1,0}} &
 \bullet  \ar@{=>}[u]_(.6){j_{0,0}}\\
 \dots \ar[r]^{i_{-2,-2}} &
 \bullet \ar[r]^{i_{-1,-2}} \ar@{=>}[u]_(.6){j_{-2,-1}} &
 \bullet \ar[r]^{i_{0,-2}} \ar@{=>}[u]_(.6){j_{-1,-1}} &
 \bullet  \ar@{=>}[u]_(.6){j_{0,-1}}\\
 \vdots\vdots\vdots &
 \vdots \ar@{=>}[u]_(.6){j_{-2,-2}} &
 \vdots \ar@{=>}[u]_(.6){j_{-1,-2}} &
 \vdots \ar@{=>}[u]_(.6){j_{0,-2}}
 }
\]

Since $\lambda_\tau$ is a $*$-rep\-re\-sen\-ta\-tion,
for each $(s,t) \in (-\bN_0)^2$, there is
a unique word $e_uf_v$ of degree $(|s|,|t|)$ such that
$\xi_{\iota_0(\mt)}$ is in the range of $\lambda_\tau(e_uf_v)$.
Set $\xi_{s,t} =  \lambda_\tau(e_uf_v)^* \xi_{\iota_0(\mt)}$.
It is not hard to see that this forms an orthonormal basis for $\H$.

For each $(s,t) \in (-\bN_0)^2$, there are unique integers
$i_{s,t} \in \{1,\dots,m\}$ and $j_{s,t} \in \{1,\dots,n\}$ so that
\begin{alignat*}{2}
 \sigma_\tau( e_{i_{s,t}}) \xi_{s-1,t} &=\xi_{s,t} &\qfor&s\le 0 \AnD t \le 0\\
 \sigma_\tau( f_{j_{s,t}}) \xi_{s,t-1} &=\xi_{s,t}  &\qfor&s \le 0 \AnD t \le 0\\
 \sigma_\tau(e_i) \xi_{s,t} &= 0 &\qif& i \ne i_{s+1,t} \OR s = 0\\
 \sigma_\tau(f_j) \xi_{s,t} &= 0 &\qif& j \ne j_{s,t+1} \OR t = 0.
\end{alignat*}
Note that we label the edges leading \textit{into} each vertex,
rather than leading out. This choice reflects the fact that
the basis vectors for a general atomic partial isometry representation
are each in the range of at most one of the partial isometries
$\pi(e_i)$ and at most one of the $\pi(f_j)$.

Consider how the  tail $\tau = e_{i_0}f_{j_0}e_{i_1}f_{j_1} \dots$
determines these integers.
It defines the path down the diagonal; that is,
\[ i_{s,s} := i_{|s|} \qand j_{s-1,s} :=j_{|s|} \qfor s \le 0 .\]
This determines the whole representation uniquely.
Indeed, for any vertex $\xi_{s,t}$ with $s,t \le 0$, take $T \ge |s|,|t|$,
and select a path from $(-T,-T)$ to $(0,0)$ that passes
through $(s,t)$.
The word $\tau_T = e_{i_0}f_{j_0} \dots e_{i_{T-1}}f_{j_{T-1}}$
satisfies $\sigma_\tau( \tau_T) \xi_{-T,-T} = \xi_{0,0}$.
Factor it as $\tau_T = w_1w_2$ with $d(w_1) = (T-|s|,T-|t|)$
and $d(w_2) = (|s|,|t|)$, so that
$\sigma_\tau(w_2) \xi_{-T,-T} = \xi_{s,t}$
and $\sigma_\tau(w_1) \xi_{s,t} = \xi_{0,0}$.
Then $w_1 = e_{i_{s,t}}w' = f_{j_{s,t}} w''$.

It is evident that each $\sigma_\tau(e_i)$
and $\sigma_\tau(f_j)$ is a partial isometry.
Moreover, each basis vector is in the range of a unique
$\sigma_\tau(e_i)$ and $\sigma_\tau(f_j)$.
So this is a defect free, atomic representation
with minimal $*$-dilation $\lambda_\tau$.
\end{exmp}

\medbreak\textbf{The symmetry group.}
An important part of the analysis of these atomic representations is the
recognition of symmetry.

\begin{defn} \label{D:shift tail}
The tail $\tau $ determines the integer data
\[ \Sigma(\tau) = \{ (i_{s,t},j_{s,t}) : s,t \le 0 \} .\]
Two tails $\tau_1$ and $\tau_2$ with data
$\Sigma(\tau_k) = \{ (i^{(k)}_{s,t},j^{(k)}_{s,t}) : s,t \le 0 \}$
are said to be \textit{tail equivalent} if the two sets of
integer data eventually coincide;
i.e.\ there is an integer $T$ so that
\[
 (i^{(1)}_{s,t},j^{(1)}_{s,t}) = (i^{(2)}_{s,t},j^{(2)}_{s,t})
 \qforal s,t \le T .
\]
Say that  $\tau_1$ and $\tau_2$ are \textit{$(p,q)$-shift tail equivalent}
for some $(p,q) \in \bZ^2$ if there is an integer $T$ so that
\[
 (i^{(1)}_{s+p,t+q},j^{(1)}_{s+p,t+q}) =
 (i^{(2)}_{s,t},j^{(2)}_{s,t}) \qforal s,t \le T .
\]
Then $\tau_1$ and $\tau_2$ are \textit{shift tail equivalent} if they
are $(p,q)$-shift tail equivalent for some $(p,q) \in \bZ^2$.

The \textit{symmetry group} of $\tau$ is the subgroup
\[
 H_\tau = \{ (p,q) \in \bZ^2 :
 \Sigma(\tau) \text{ is $(p,q)$-shift tail equivalent to itself} \} .
\]
A sequence $\tau$ is called \textit{aperiodic} if $H_\tau = \{(0,0)\}$.

The semigroup $\Fth$ is said to satisfy the
\textit{aperiodicity condition} if there is an aperiodic infinite word.
Otherwise we say that $\Fth$ is \textit{periodic}.
\end{defn}

In our classification of atomic $*$-rep\-re\-sen\-ta\-tions,
an important step is to define a symmetry group
for the more general representations which occur.
The $*$-rep\-re\-sen\-ta\-tion will be irreducible precisely when
the symmetry group is trivial.
This will yield a method to decompose atomic
$*$-rep\-re\-sen\-ta\-tions as direct integrals of
irreducible atomic $*$-rep\-re\-sen\-ta\-tions.

\medbreak\textbf{The graph of an atomic representation.}
We need to develop a bit more notation.
Let $\sigma$ be an atomic representation.
Let the corresponding basis be $\{\xi_k : k \in S\}$.
Write $\dot\xi_k$ to denote the subspace $\bC \xi_k$.
Form a graph $\G_\sigma$ with vertices $\dot\xi_k$.
If $\sigma(e_i) \dot\xi_k = \dot\xi_l$, draw a directed \textit{blue} edge
from $\dot\xi_k$ to $\dot\xi_l$ labelled $i$; and if
$\sigma(f_j) \dot\xi_k = \dot\xi_l$, draw a directed \textit{red} edge
from $\dot\xi_k$ to $\dot\xi_l$ labelled $j$.
This is the graph of the representation, and it contains all
of the information about $\sigma$ except for the scalars in $\bT$.

In our analysis of atomic representations, one can easily split
a representation into a direct sum of atomic representations
which have a connected graph. So we will generally work with
representations with connected graph.

\begin{lem}\label{L:push pull}
Let $\sigma$ be a defect free atomic representation
with connected graph $\G_\tau$.
Let $\dot\xi_1$ and $\dot\xi_2$ be two vertices in $\G_\sigma$.
Then there is a vertex $\dot\eta$ and words $w_1,w_2 \in \Fth$ so that
$\dot\xi_i = \sigma(w_i) \dot\eta$ for $i=1,2$.
\end{lem}

\begin{pf}
The connectedness of the graph means that there is a path
from $\dot\xi_1$ to $\dot\xi_2$.
We will modify this path to first pull back along a path leading into $\dot\xi_1$,
and then move forward to $\dot\xi_2$.

The original path can be written formally as
$a_k a_{k-1} \dots a_1$ where each $a_l$ has the form $e_i$ or $f_j$
if it is moving forward or $e_i^*$ or $f_j^*$ if pulling back.
After deleting redundancies, we may assume that there are no adjacent
terms $f_j^* f_j$ or $e_i^* e_i$.
At each vertex, there is a unique blue (red) edge leading in;
so moving forward along a blue (red) edge and pulling back
along the same colour is just one of these redundancies.
Thus if there are adjacent terms of the form $a^*b$ which do not cancel,
then one of $a,b$ is an $e$ and the other an $f$.

For definiteness, suppose that this section of the path moves
from $\dot\eta_1$ to $\dot\eta_2$ along the path $f_j^*e_i$.
That means that
$\sigma(e_i)\dot \eta_1 = \dot\eta = \sigma(f_j)\dot\eta_2$.
As $\sigma$ is defect free, there is a basis vector $\dot\eta_0$
and an $f_{j'}$ so that $\sigma(f_{j'})\dot\eta_0 = \dot\eta_1$.
Factor $e_i f_{j'}$ in the other order as $f_{j''}e_{i'}$.
Then
\[
 \sigma(f_j)\dot\eta_2 = \dot\eta = \sigma(e_i f_{j'}) \dot\eta_0
 = \sigma(f_{j''}e_{i'}) \dot\eta_0 = \sigma(f_{j''}) \big( \sigma(e_{i'}) \dot\eta_0 \big) .
\]
However there is a unique red edge into $\eta$, and thus
\[ f_{j''} = f_j \qand \sigma(e_{i'}) \dot\eta_0 = \dot\eta_2 .\]
It is now clear that there is a red--blue diamond in the graph
with apex $\dot\eta_0$, and red and blue edges leading to
$\dot\eta_1$ and $\dot\eta_2$ respectively, and from there,
red and blue edges respectively leading into $\dot\eta$.
So the path $f_j^*e_i$ may be replaced by $e_{i'}f_{j'}^*$.
Similarly, the path $e_i^*f_j$ from $\dot\eta_2$ to $\dot\eta_1$
may be replaced by $f_{j'}e_{i'}^*$.

Repeated use of this procedure replaces any path from $\dot\xi_1$ to $\dot\xi_2$
by a path of the form $w_2w_1^*$; and hence $\eta = \sigma(w_1)^*\xi_1$
is the intermediary vector.
\bx\end{pf}

\section{Classifying Atomic Representations} \label{S:atomic}

Consider an atomic $*$-rep\-re\-sen\-ta\-tion $\pi$ of $\Fth$
with connected graph $\G_\pi$. Observe that if we restrict
$\pi$ to the subalgebra generated by the $e_i$'s, we obtain an
atomic, defect free representation of the free semigroup $\bF_m^+$.
The graph splits into the union of its blue components. By
\cite{DP1}, this decomposes the restriction into a direct sum of
ring representations and infinite tail representations. Our first
result shows that the connections provided by the red edges force a
parallel structure in the blue components.

\begin{lem}\label{L:parallel}
Let $\pi$ be an atomic $*$-rep\-re\-sen\-ta\-tion of $\Fth$
with connected graph $\G_\pi$. Let $H_1$ and $H_2$ be two blue
components of $\G_\pi$ and assume that there is a red edge
leading from a vertex in $H_1$ to a vertex in $H_2$.
Then every red edge into $H_2$ comes from $H_1$.
The two components  are either both of infinite tail type
or both are of ring type; and in the ring case the length
of the (unique) ring for $H_2$ is an integer multiple $t$ of
the ring length in $H_1$, where $1 \le t \le n$.
\end{lem}

\begin{pf}
This is an exercise in using the commutation relations.
Suppose first that $H_1$ is an infinite tail graph.
Fix a vertex $\dot\xi_0 \in H_1$ and a red edge $f_{j_0}$ so that
$\pi(f_{j_0}) \dot\xi_0 = \dot\zeta_0$ is a vertex in $H_2$.
There is a unique infinite sequence of vertices $\xi_k$ for $k<0$
and integers $i_k$ so that
\[ \pi(e_{i_k})\dot\xi_{k-1} = \dot\xi_k \qfor k \le 0 .\]
Now there are unique integers $i'_k$ and $j_k$ so that
\[ f_{j_0} e_{i_0i_{-1}\dots i_{k+1}} = e_{i'_0i'_{-1}\dots i'_{k+1}} f_{j_k} \qfor k < 0 .\]
Let $\dot\zeta_k := \pi( f_{j_k}) \dot\xi_k$ for $k<0$.
Then it is evident that
\[ \pi(e_{i'_k})\dot\zeta_{k-1} = \dot\zeta_k \qfor k \le 0 .\]
The images of each vertex under the various red edges are all distinct.
So in particular, the $\dot\zeta_k$ are all distinct,
and so $H_2$ is also an infinite tail graph.

Now one similarly can follow each vertex $\dot\zeta_k$ forward under
a blue path $e_u$ to reach any vertex $\dot\zeta$ in $H_2$.
Then
\[
 \dot\zeta = \pi(e_u)\dot\zeta_k
 = \pi(e_u) \pi( f_{j_k}) \dot\xi_k
 = \pi(f_{j'}) \pi(e_{u'})  \dot\xi_k.
\]
Thus the vertex $\pi(e_{u'})  \dot\xi_k$ in $H_1$ is mapped
to $\dot\zeta$ by $\pi(f_{j'})$.
Hence the red edge leading into $\dot\zeta$ comes from $H_1$.

The case of a ring representation is similar.
Starting with any vertex in $H_1$ which is connected to $H_2$ by a red edge,
one can pull back along the blue edges until one is in the ring.
So we may suppose that $\dot\xi_0$ lies in the ring of $H_1$
and that $\pi(f_{j_0}) \dot\xi_0 = \dot\zeta_0$ in $H_2$.
Let $u$ be the unique minimal word such that
$\pi(e_u) \dot\xi_0 = \dot\xi_0$, and let $p=|u|$.

As in the first paragraph, we continue to pull back from $\dot\zeta_0$
and from $\dot\xi_0$ along the blue edges.
Call these edges $\dot\zeta_k$ and $\dot\xi_k$ respectively for $k<0$.
After $p$ steps, we return to $\dot\xi_0 = \dot\xi_{-p}$ in $H_1$,
and we have obtained $p$ distinct vertices in $H_2$ and
reach $\dot\zeta_{-p}$.
So there is a word $u'$ of length $p$ so that
$\pi(e_{u'})\dot\zeta_{-p} = \dot\zeta_0$.
The commutation relations yield
\[
 \pi(e_{u'})\dot\zeta_{-p} = \dot\zeta_0 = \pi(f_{j_0}) \pi(e_u) \dot\xi_0
 = \pi(e_{u''}) \pi(f_{j'}) \dot\xi_0 .
\]
There is a unique blue path of length $p$ leading into $\dot\zeta_0$.
Therefore $u''=u'$ and $\pi(f_{j'}) \dot\xi_0 = \dot\zeta_{-p}$.
Notice that if $j'=j_0$, then $\dot\zeta_{-p} = \dot\zeta_0$.
However if $j'\ne j_0$, then $\dot\zeta_{-p}$ is a different vertex in $H_2$.

Repeat the process, pulling back another $p$ blue steps, to reach
a vertex $\dot\zeta_{-2p}$.
Along the way, we obtain vertices which are distinct from the previous
ones, and each is the image of some vertex in the ring of $H_1$.
As before, $\dot\zeta_{-2p}$ is the image of $\dot\xi_0$ under some
red edge.
Eventually this process must repeat, because there are only $n$
red edges out of $\dot\xi_0$.
That is, there are integers $s$ and $t$ with $1 \le t \le n$ so that
$\dot\zeta_{-(s+t)p} = \dot\zeta_{-sp}$.
This is a ring in $H_2$ of length $tp$.

The argument that each edge in $H_2$ is in the range of a red edge
coming from $H_1$ is identical to the infinite tail case.
\bx\end{pf}

One might hope that the red edges from $H_1$ to $H_2$ provide a nice
bijection, or a $t$-to-$1$ map that preserves the graph structure.
Even though these edges are determined algebraically by the
commutation relations, such a nice pairing does not occur
as the following examples demonstrate.

\begin{exmp}
Consider $m=n=3$ with $\theta$ given by $\big( (1,2) , (2,1) \big)$,
or equivalently by the relations
\[
 e_1f_2 = f_1e_2, \quad e_2f_1=f_2e_1,
 \qand e_if_j = f_j e_i \quad\text{otherwise.}
\]
There is a 1-dimensional defect free representation $\rho(e_3) = \rho(f_3) = 1$
and $\rho(e_i) = \rho(f_i) = 0$ for $i=1,2$.
This has a dilation to a $*$-repre\-sent\-ation $\pi$.
Let the initial vector be called $\xi_0$.
Define $\zeta_{j0} = \pi(f_j)\xi_0$ for $j=1,2$;
and let $\xi_i = \pi(e_i) \xi_0$ and $\zeta_{ji} = \pi(e_i)\zeta_{j0}$
for $i=1,2$ and $j=1,2$.
The blue component $H_0$ containing $\xi_0$ has a ring of length one.
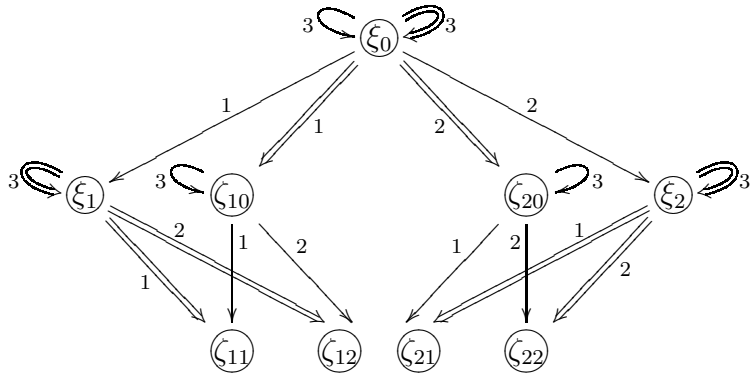
\[
 \xymatrix@!@R=.1in@C=0.05in{
&&{\xy *++={\xi_0} *\frm{o} \endxy}
 \ar@(ul,l)[]_(.6){3} \ar@{=>}@(ur,r)[]^(.6){3}
\ar[dll]_-1 \ar@{=>}[dl]^-1 \ar@{=>}[dr]_-2 \ar[drr]^-{2}&&
 \\
{\xy *++={\xi_1} *\frm{o} \endxy}
\ar@{=>}@(ul,l)[]_(.6)3 \ar@{=>}[dr]_-1  \ar@{=>}[drr]^(.3){2}
& {\xy *++={\zeta_{10}} *\frm{o} \endxy}
\ar@(ul,l)[]_(.6)3 \ar[d]^(.3)1 \ar[dr]^(.4)2
&&{\xy *++={\zeta_{20}} *\frm{o} \endxy}
\ar@(ur,r)[]^(.6)3 \ar[dl]_(.4)1 \ar[d]_(.3)2
& {\xy *++={\xi_2} *\frm{o} \endxy}
\ar@{=>}@(ur,r)[]^(.6)3 \ar@{=>}[dl]^(.4)2 \ar@{=>}[dll]_(.3)1
\\
&{\xy *++={\zeta_{11}} *\frm{o} \endxy}
&{\xy *++={\zeta_{12}} *\frm{o} \endxy}
\quad
{\xy *++={\zeta_{21}} *\frm{o} \endxy}
&{\xy *++={\zeta_{22}} *\frm{o} \endxy}
}
\]

It is easy to check that $\pi(e_3)\zeta_{j0} = \zeta_{j0}$ for $j=1,2$.
Hence each is a ring in a separate blue component $H_j$.
But the commutation relations also show that $\pi(f_j) \xi_i = \zeta_{ij}$.
This means that the vertex $\dot\xi_1$ has two red edges leading to the
component $H_1$; and $\dot\xi_2$ has two red edges leading to $H_2$.
\end{exmp}

\begin{exmp}
With the same algebra, consider the 1-dimensional representation
$\rho(e_1) = \rho(f_3) = 1$ and $\rho(e_i) = \rho(f_j) = 0$ otherwise.
Again this has a dilation to an isometric defect free representation $\pi$.
Let the initial vector be called $\xi$.
Define $\zeta_j = \pi(f_j) \xi$ for $j=1,2$.
A computation with the relations shows that $\pi(e_1) \zeta_1 = \zeta_1$.
So this is a ring of length one in a component $H_1$.
But $\pi(e_2) \zeta_1 = \zeta_2$.
So $\zeta_2$ is also in $H_1$; but even though it is in the range
of a red edge from the ring of $H_0$, it does not lie on the ring of $H_1$.
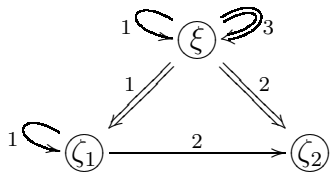
\[
 \xymatrix{&{\xy *++={\xi} *\frm{o} \endxy}
  \ar@(ul,l)[]_(.6){1} \ar@{=>}@(ur,r)[]^(.6){3}
  \ar@{=>} [dl]_-1 \ar@{=>}[dr]^-{2}
 \\
 {\xy *++={\zeta_1} *\frm{o} \endxy} \ar[rr]^-{2}  \ar@(ul,l)[]_(.6){1} &
 &{\xy *++={\zeta_2} *\frm{o} \endxy} }
\]
\end{exmp}

\begin{exmp}\label{E:reverse favourite}
Consider $m=n=2$ and the reverse 3-cycle semigroup
of Example~\ref{E:favourite} given by the 3-cycle
$\big( (1,1) , (2,1) , (1,2) \big)$.
There is a 1-dimensional defect free representation
\[ \rho(e_2) = \rho(f_2) = 1 \qand  \rho(e_1) = \rho(f_1) = 0 .\]
This has a dilation to a $*$-rep\-re\-sen\-ta\-tion $\pi$.
Let the initial vector be called $\xi$.
Define $\eta = \pi(f_1)\xi$ and $\zeta_j =  \pi(f_j)\eta$ for $j=1,2$.
Then an exercise with the relations shows that
$\pi(e_1)\eta = \eta$ and $\pi(e_1)\zeta_2 = \zeta_1$
and $\pi(e_2)\zeta_1 = \zeta_2$.
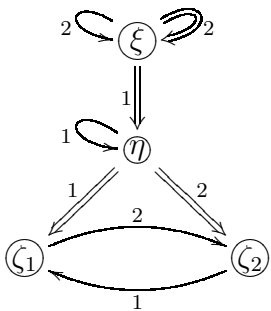
\[
 \xymatrix{
&{\xy *++={\xi} *\frm{o} \endxy}
 \ar@(ul,l)[]_(.6)2 \ar@{=>}@(ur,r)[]^(.6)2 \ar@{=>}[d]_-1 &
\\ &{\xy *++={\eta} *\frm{o} \endxy}
\ar@(ul,l)[]_(.6)1 \ar@{=>}[dl]_-1 \ar@{=>}[dr]^-2
\\ {\xy *++={\zeta_1} *\frm{o} \endxy}
\ar@/^{1pc}/[rr]^-2 && {\xy *++={\zeta_2} *\frm{o} \endxy}
\ar@/^{1pc}/[ll]^-1
}
\]
Thus the initial component $H_0$ has a ring of length one at $\xi$,
as does the component $H_1$ connected to it (at $\eta$), but it
connects to a component $H_2$ in which $\zeta_i$ are the vertices
of a ring of length 2.
One can show that there are components with rings of length $2^k$
for all $k\ge0$.
\end{exmp}

\medbreak\textbf{Splitting into cases.}
Evidently the reasoning of Lemma~\ref{L:parallel} also applies when
we decompose the graph into its red components.
It is now possible to split the analysis into several cases.

\textbf{1. ring by ring type.}
\textit{Both blue and red components are ring representations. }
In this case, the set of vertices which are in both a red and a blue ring
determines a finite dimensional coinvariant subspace on which the
representation is defect free, and there is exactly one red and blue edge
beginning at each vertex.
Moreover this is a cyclic subspace for the representation because
starting at any basis vector, pulling back along the blue and red
edges eventually ends in the ring by ring portion.
So this finite dimensional piece determines the full representation.

Since each ring is obtained by pulling back from any of the others,
it follows from Lemma~\ref{L:parallel} that all of the blue rings have
the same length, say $k$;
and likewise all of the red rings have the same length, say $l$.

\textbf{2. mixing type.}
\textit{For one colour, the components are ring representations while
the for the other colour, the components are infinite tail representations.}
There are two subcases.

\textbf{2a. ring by tail type.}
\textit{The blue components are ring representations, and the red components
are infinite tail type.}

If one begins at any blue ring and pulls back along the red edges,
one obtains an infinite sequence of blue ring components.
By Lemma~\ref{L:parallel}, the size of the ring is decreasing as
one pulls back.  Hence it is eventually constant.
 From this point on back, there is a unique red edge from each ring
to the corresponding point on the next.
Thus one obtains a semi-infinite cylinder of fixed circumference $k$
which is coinvariant and cyclic, and  determines the full representation.

\textbf{2b. tail by ring type.}
\textit{The red components are ring representations, and the blue components
are infinite tail type.}

\textbf{3.  tail by tail type.}
\textit{Both the red and the blue components are of infinite tail type.}
This actually has several subtypes.

Start with a basis vector $\dot\xi$ in a blue component $H_0$.
Pull back along the red edges to get an infinite sequence of blue
components $H_t$ for $t\le0$.
The union of these components for $t \le 0$ forms a coinvariant cyclic
subspace that determines the full representation.
This sequence may be eventually periodic, or they may all be distinct.
If they are eventually periodic, we may assume that we begin with a
component $H_0$ in the periodic sequence.

\textbf{3a. inductive type.}
\textit{The sequence of components are all distinct.}

In this case, starting at any basis vector $\dot\xi_{0,0}$, one
may pull back along both blue and red edges to obtain basis
vectors $\xi_{s,t}$ for $(s,t) \in (-\bN_0)^2$.
The restriction of the representation to this coinvariant subspace
is defect free, and determines the whole representation
as in Example~\ref{3a_repn}. So this is an inductive representation.

\textbf{3b.} \textit{The sequence of components repeats after $l$ steps. }
Thus by Lemma~\ref{L:parallel}, there are blue components
$H_0, \dots,H_{l-1}$ so that every red edge into each $H_i$ comes from
$H_{i-1 \!\!\pmod l}$.
This is further refined by comparing the point of return
to the starting point $\dot\xi$.
It is not apparent at this point that the only possibilities are the following:

\textbf{3bi. return below.}
\textit{The return is eventually below the start.}  In this case,
there is a vertex $\dot\xi_0$ in $H_0$ and a word $u_0$
so that the vertex $\dot\zeta_0$ obtained by pulling back $l$ red steps
using the word $v_0$ from $\dot\xi_0$ satisfies
\[ \pi(e_{u_0}) \dot\xi_0 = \dot\zeta_0 \qand \pi(f_{v_0}) \dot\zeta_0 = \dot\xi_0. \]
This same type of relationship persists when pulling back along both blue
and red edges from $\dot\xi_0$.

\textbf{3bii. return above.}
\textit{The return is eventually above the start.}
Here there is a vertex $\dot\xi_0$ in $H_0$ and a word $u_0$
so that the vertex $\dot\zeta_0$ obtained by pulling back $l$ red steps
using the word $v_0$ from $\dot\xi_0$ satisfies
\[ \pi(e_{u_0}) \dot\zeta_0 = \dot\xi_0 = \pi(f_{v_0}) \dot\zeta_0 . \]
This same type of relationship persists when pulling back along both blue
and red edges from $\dot\xi_0$.
\smallskip

The ring by ring representations are \textit{finitely correlated},
meaning that there is a finite dimensional coinvariant, cyclic subspace.
Equivalently, this means that the $*$-rep\-re\-sen\-ta\-tion is the minimal
$*$-dilation of a finite dimensional defect free representation.

Conversely, every finitely correlated atomic  $*$-rep\-re\-sen\-ta\-tion
with connected graph is of ring by ring type, as the other cases
clearly do not have a finite dimensional
non-zero coinvariant and cyclic subspace.

In the following sections, each case will be considered in more detail.
Eventually a common structure emerges.
This will be codified by the general construction given in the next section.

\medbreak\textbf{Symmetry.}
In each case, we associate a symmetry group to the picture.
In the ring by ring case, it is easiest to describe because
there is no equivalence relation.
If the minimal coinvariant subspace consists of
blue cycles of length $k$ and red cycles of length $l$,
then we associate the representation to a quotient
group $G$ of $\C_k \times \C_l$;
and the symmetry group is a subgroup of $G$.
We show that there is a finite dimensional representation
on $\bC^{kl}$ which reflects the full symmetry, and that
certain quotients yield a decomposition into
irreducible summands.

In type 3a, the inductive case, we have already seen how
to define a symmetry subgroup of $\bZ^2$.
Again, an inductive representation will be irreducible
precisely when this symmetry group is trivial.

In the other cases, the symmetry group is a
subgroup of $\C_k \times \bZ$ in the type 2a case
or of $\bZ^2/\bZ(a,b)$ in the 3b cases.

In type 2 and 3, the decomposition into irreducibles may
require a direct integral rather than a direct sum.

\section{A group construction}
\label{S:group}

In this section, we will describe a general class of examples,
and explain how to decompose them into irreducible representations.

Start with an abelian group $G$ with two designated generators
$\fg_1$ and $\fg_2$.
We consider a defect free atomic representation of $\Fth$ on $\ltwo(G)$
which is given by the following data:
\begin{alignat*}{2}
 &i:G\to\{1,\dots,m\} &\qquad i(g) &=: i_g\\
 &j:G\to\{1,\dots,n\} &\qquad j(g) &=: j_g\\
 &\alpha:G\to\bT &\qquad \alpha(g) &=: \alpha_g\\
 &\beta:G\to\bT &\qquad \beta(g) &=: \beta_g
\end{alignat*}
We wish to define a representation $\sigma$ of $\Fth$ by
\begin{align*}
 \sigma(e_i) \xi_g &= \delta_{i\,i_g}\, \alpha_g \xi_{g+\fg_1}\\
 \sigma(f_j) \xi_g &= \delta_{j\,j_g}\, \beta_g \xi_{g+\fg_2} .
\end{align*}
In order for this to actually be a representation, we require that
\[ e_{i_{g+\fg_2}} f_{j_g} = f_{j_{g+\fg_1}}e_{i_g} \qforal g \in G \]
and
\[
 \alpha_{g+\fg_2} \beta_{g} = \beta_{g+\fg_1} \alpha_{g}
 \qforal g \in G.
\]
Let $\G_\sigma$ be the graph of this representation, which has
vertices $\dot\xi_g$ for $g \in G$ and blue edges labelled $i_g$
from $\dot\xi_g$ to $\dot\xi_{g+\fg_1}$, and red edges labelled
$j_g$ from $\dot\xi_g$ to $\dot\xi_{g+\fg_2}$.
These will be called \textit{group construction representations} of $\Fth$.

It is evident that there is a unique blue and red edge leading into
each vertex $\dot\xi_g$, and so this is a defect free atomic representation.
Thus it has a unique minimal $*$-dilation.
Decomposing this representation into irreducible summands
or a direct integral must simultaneously decompose the
$*$-dilation into a direct sum or direct integral of the
minimal $*$-dilations of the summands.
The symmetry group of $\sigma$ is defined as
\[
 H = \{h \in G : i_{g+h} = i_{g}
 \AnD j_{g+h} = j_{g} \FORAL g \in G \} .
\]
This will play a central role in this decomposition.

We will address the non-trivial issue of how to actually define
the functions $i$ and $j$ later when considering the various cases.
There are obstructions, and the way by which this difficulty
is overcome is not immediately apparent.
Example~\ref{E:cycle reps} is an example of this type in which $G$
is the cyclic group $\C_k$ and $\fg_1 = -\fg_2= 1$.
The case of $G=\bZ^2$ can be seen in Example~\ref{3a_repn}.
Here we describe a coinvariant subspace which is identified
with $(-\bN_0)^2$, but it can be extended to all of $\bZ^2$.
(We apologize to the reader that the notation is not consistent
between Example~\ref{3a_repn} and this section.)

The issue of the scalar functions $\alpha$ and $\beta$ is more
elementary.  We shall see that it suffices to consider the case
in which $\alpha$ and $\beta$ are constant, and will determine
when two are unitarily equivalent.
For the moment, we assume that such a representation is given
and consider how to analyze it.

The group $G$, being abelian with two generators, is a
quotient of $\bZ^2$.  The subgroups of $\bZ^2$ are $\{0\}$,
singly generated $\bZ(a,b)$, or doubly generated.  In the
doubly generated case, it is an easy exercise to see that
the quotient is finite.  In this case, if $\fg_i$ has order $p_i$
for $i=1,2$, then $G$ is in fact a quotient of $\C_{p_1}\times \C_{p_2}$.

We shall show that the different groups correspond to
the different representation types as follows:
\newcounter{type}
\renewcommand{\thetype}{Type \arabic{type}. }
\begin{list}{\thetype\hfill}
{\usecounter{type} \rightmargin=0pt\leftmargin=72pt%
 \labelwidth=46pt\labelsep=2pt\itemsep=1pt}
\item \qquad $G$ finite.
\item (a)\ \, $G = \C_k \times \bZ$;\\(b)\ \,  $G = \bZ \times \C_l$.
\item (a)\ \, $G=\bZ^2$;\\
(bi)\, $G=\bZ^2/\bZ(k,l)$, $kl>0$;\\
(bii) $G=\bZ^2/\bZ(k,l)$, $kl<0$.
\end{list}

\medbreak\textbf{Scalars.}
Let us dispense with the scalars first.
This is actually quite straightforward in spite
of the notation.
For each group $G$, there is a canonical homomorphism $\kappa$
of $\bZ^2$ onto $G$ sending the standard generators
$(1,0)$ and $(0,1)$ to $\fg_1$ and $\fg_2$.
Let $K$ be the kernel of $\kappa$, so that $G \simeq \bZ^2/K$.

\begin{thm} \label{T:scalars}
Let $G = \bZ^2/K$ as above.
A group construction representation $\sigma$ on $\ltwo(G)$
is unitarily equivalent to another of the
same type with the same functions $i,j$ for which
the scalar functions are constants  $\alpha_0$ and $\beta_0$.

The constants determine a unique character $\psi$ of $K$.
They are unique up to a choice of an extension $\phi$ of $\psi$
to a character on $\bZ^2$; and they are
given by $\alpha_0 = \phi(1,0)$ and $\beta_0 = \phi(0,1)$.
The choice of $\phi$ is unique up to a character $\chi$ of $G$,
which changes $\alpha_0$ and $\beta_0$
to $\alpha_0\chi(\fg_1)$ and $\beta_0\chi(\fg_2)$.
\end{thm}

\begin{pf}
For each $(s,t) \in \bZ^2$, there is a unique path from $\dot\xi_0$
to $\dot\xi_{\kappa(s,t)}$ in $\G_\sigma$ with $s$ blue edges
and $t$ red edges, where a negative number indicates traversing the
arrow in the backward direction.  Corresponding to this, there is a
unique partial isometry $W$ in $\ca(\sigma(\Fth))$ of the form
$\sigma(e_u)\sigma(f_v)$, $\sigma(f_v)^*\sigma(e_u)$,
$\sigma(e_u)^*\sigma(f_v)$ or $\sigma(e_u)^*\sigma(f_v)^*$
with $|u|=|s|$ and $|v|=|t|$, depending on the signs of $s$ and $t$,
so that $W \dot \xi_0 = \dot \xi_{\kappa(s,t)}$.
The subgroup $K$ corresponds to those paths which return
$\dot\xi_0$ to itself.
Thus for $(s,t) \in K$,  $W \xi_0 = \psi(s,t) \xi_0$  for a unique scalar
$\psi(s,t) \in \bT$.

In any unitarily equivalent representation given in the same form
on $\ltwo(G)$, the vector sent to $\xi_0$ must be a vector with
exactly the same functions $i,j$; and thus must come from
$\spn\{ \xi_h : h \in H\}$.  Moreover there are constraints on
the scalars, namely if $U^*\xi_0 = \sum_{h \in A} a_h \xi_h$
where $A = \{h : a_h \ne 0\}$, then $\alpha(h+g) = \alpha(h'+g)$
and $\beta(h+g) = \beta(h'+g)$ for all $h,h' \in A$ and all $g \in G$.
(For otherwise, the representation would not correspond to a graph.)
Hence exactly the same words return each $\xi_h$ to itself, for $h \in A$.

We claim that the function $\psi$ on $K$ is independent of
the unitary equivalence.
Let $\psi'$ be the function obtained in this equivalent representation.
The argument of the previous paragraph shows that instead
of computing the function $\psi'$ using $U^*\xi_0$, we can
use $\xi_h$ for any $h \in A$ and obtain the same result.

Let $W$ be the partial isometry found in the first paragraph which
carries $\dot\xi_0$ to $\dot\xi_h$.
Let $(s,t) \in K$, and let $W_0$ and $W_h$ be the partial isometries
of degree $(s,t)$ which take $\dot\xi_0$ and $\dot\xi_h$ to
themselves, respectively.
Now $W^*WW_0$ and $W^*W_hW$ are partial isometries of the
same combined and absolute degrees which map $\dot\xi_0$
to itself.
The uniqueness of factorization means that these two words
are equal!
Therefore
\[
 \psi(s,t) = \ip{W^*WW_0 \xi_0,\xi_0}
 = \ip{W^*W_hW \xi_0,\xi_0} = \psi'(s,t) .
\]
So we see that the original choice of scalars forces certain
values, namely the function $\psi$, to be fixed independent
of unitary equivalence.

Next we show that $\psi$ is a character.
Given two words $(s_1,t_1)$ and $(s_2,t_2)$ in $K$,
let $W_1$ and $W_2$ be the corresponding partial isometries.
Then the partial isometry $W$ corresponding to the sum
\mbox{$(s_1 \!+\! s_2,t_1 \!+\! t_2)$}
need not equal $W_1W_2$.
However the unique factorization means that
$W_1W_2 =V^*VW$ for some partial isometry $V$
containing $\xi_0$ in its domain.
Thus one computes
\begin{align*}
 \psi(s_1+s_2,t_1+t_2) &= \ip{V^*VW \xi_0,\xi_0} \\
 &=\ip{W_1W_2 \xi_0,\xi_0} =\psi(s_1,t_1)\psi(s_2,t_2) .
\end{align*}
So $\psi$ is multiplicative.

It is routine to extend $\psi$ to a character $\phi$ of $\bZ^2$.
In our context, one can easily do this `bare hands'.
But it is a general fact for characters on any subgroup
of any abelian group \cite[Corollary~24.12]{HR}.
If $\phi_1$ and $\phi_2$ are two characters extending $\psi$,
then $\phi_2\ol{\phi_1}$ is a character which takes the constant
value 1 on all of $K$.  Thus it induces a character $\chi$ of $G$.
We see that $\phi_2(s,t) = \phi_1(s,t)\chi(s\fg_1+t\fg_2)$.
Conversely, any choice of $\chi$ yields an extension.

The unitary equivalence between $\sigma$ and the representation
$\sigma'$ on the same graph but with constants $\alpha_0=\phi(1,0)$
and $\beta_0=\phi(0,1)$ can be accomplished by a diagonal
unitary $U = \diag(\gamma_g)$.  Define $\gamma_g$
by selecting a partial isometry $W$ as in the first
paragraph so that $W\dot\xi_0 = \dot\xi_g$.
Let $(a,b)$ be the degree of the word $W$.
Define
\[ \gamma_g = \phi(a,b) \ol{\ip{W\xi_0,\xi_g}} .\]
While the choice of $W$ is not unique, if $W'$ is another such
partial isometry, $W^*W'\dot\xi_0 = \dot\xi_0$ corresponds
to a word $(s,t)$ in $K$.  Then $W'$ has degree $(a+s,b+t)$.
Thus
\begin{align*}
 \phi(a+s,b+t) \ol{\ip{W'\xi_0,\xi_g}} &=
 \phi(a,b) \phi(s,t) \ol{\ip{W W^*W'\xi_0,\xi_g}} \\ &=
 \phi(a,b) \phi(s,t)\ol{\phi(s,t)} \ol{\ip{W\xi_0,\xi_g}}\\ &=
 \phi(a,b) \ol{\ip{W\xi_0,\xi_g}}
\end{align*}
So $U$ is well defined.

Now if $\gamma_g$ is computed using $W$, calculate
\begin{align*}
 \gamma_{g+\fg_1} &=
 \phi(a+1,b) \ol{\ip{\sigma(e_{i_g})W\xi_0,\xi_{g+\fg_1}}} \\ &=
 \phi(a+1,b) \ol{\alpha_g} \ol{\ip{W\xi_0,\xi_g}} .
\end{align*}
Then we see that
\begin{align*}
 U\sigma(e_{i_g})U^* \xi_g &=
 U\sigma(e_{i_g}) \ol{\phi(a,b)} \ip{W\xi_0,\xi_g} \xi_g \\ &=
 U \alpha_g \ol{\phi(a,b)} \ip{W\xi_0,\xi_g} \xi_{g+\fg_1}\\ &=
 \phi(a+1,b) \ol{\alpha_g} \ol{\ip{W\xi_0,\xi_g}} \alpha_g
 \ol{\phi(a,b)} \ip{W\xi_0,\xi_g} \xi_{g+\fg_1}\\ &=
 \phi(1,0) \xi_{g+\fg_1} = \alpha_0 \xi_{g+\fg_1} .
\end{align*}
One deals with $\sigma(f_{j_g})\xi_g$ in the same manner.

Thus $\sigma$ is unitarily equivalent to the representation
$\sigma'$ with scalars $\alpha_0$ and $\beta_0$.
It was irrelevant which extension $\phi$ of $\psi$ was
used; so all are unitarily equivalent to each other.
\bx\end{pf}

\medbreak\textbf{Symmetry.}
The key to the decomposition is to look for symmetry
in the graph of $\sigma$.  Recall that
\[
 H = \{h \in G : i_{g+h} = i_{g}
 \AnD j_{g+h} = j_{g} \FORAL g \in G \} .
\]
It is clear that $H$ is a subgroup of $G$.
Note that we ignore the scalars for this purpose
which is justified by Theorem~\ref{T:scalars}.
Indeed, we shall suppose that the scalars are constants
$\alpha_0$ and $\beta_0$.

For each coset $[g] = g+H$ of $G/H$, let
\[ \W_{[g]} = \spn\{ \xi_{g+h} : h \in H \} .\]
Pick a representative $g_k$ for each coset of $G/H$,
selecting $0 \in [0]$.
Each of the subspaces $\W_{[g_k]}$ can be identified with $\ltwo(H)$.
This identification depends on the choice of representative.
Define $J_{[g_k]}:\ltwo(H) \to \W_{[g_k]}$ by $J_{[g_k]} \xi_h = \xi_{g_k+h}$.

Since $i_g=i_{g+h}$ and $j_g=j_{g+h}$, these integers are
dependent only on the coset.  So we will write $i_{[g]}$ and $j_{[g]}$.
Observe that $\sigma(e_{i_{[g]}})$ maps $\W_{[g]}$ to $\W_{[g+\fg_1]}$.
To understand this map, note that for each $g_k$, there is a unique element
$h_{k1} \in H$ so that the representative for
$[g_k+\fg_1]$ is $g_k+\fg_1 - h_{k1}$.
Then it is easy to see that
\[ \sigma(e_{i_{[g_k]}}) |_{\W_{[g_k]}} = \alpha_0 J_{[g_k+\fg_1]} L_{h_{k1}} J_{[g_k]}^* \]
where $L_h$ is the (left) regular translation by $h$ on $\ltwo(H)$.

Likewise, there is an $h_{k2} \in H$ so that $g_k+\fg_2 - h_{k2}$
is the chosen representative for $[g_k+\fg_2]$.  Then
\[ \sigma(f_{j_{[g_k]}}) |_{\W_{[g_k]}} = \beta_0 J_{[g_k+\fg_2]} L_{h_{k2}} J_{[g_k]}^* .\]

It is now a routine matter to diagonalize $\sigma$.
The unitary operators $L_{h}$ for $h \in H$ all commute, and so
can be simultaneously diagonalized  by the Fourier transform
which identifies $\ltwo(H)$ with $L^2(\hat{H})$, and carries
$L_h$ to the multiplication operator $M_h$ given by
$M_h f(\chi) = \chi(h) f(\chi)$.
\smallbreak

We explain in more detail how this works in the finite case.
Here $L_h$ are unitary matrices, and $\sigma$ decomposes as
a finite direct sum of irreducible representations.
For each $\chi \in \hat{H}$, let
\[ \zeta_0^\chi = |H|^{-1/2} \sum_{h\in H} \ol{\chi(h)} \xi_h .\]
Then a routine calculation shows that $L_h \zeta_0^\chi = \chi(h) \zeta_0^\chi$.

Consider the subspaces
\[
 \M_\chi = \spn \{ \zeta_{[g]}^\chi := J_{[g]} \zeta_0^\chi
 : [g] \in G/H \} .
\]
The choice of representative for each coset only affects the scalar multiple of the
vectors, and the subspace $\M_\chi$ is independent of this choice.
It is easy to see that these are reducing subspaces for $\sigma(\Fth)$.
This decomposes $\sigma$ into a direct sum of representations
$\sigma_\chi$ acting on $\ltwo(G/H)$.
Indeed, if $[g_k+\fg_1]=[g_l]$, then we calculate
\begin{align*}
 \sigma_\chi(e_{i_{[g_k]}}) \zeta_{[g_k]}^\chi &=
 \sigma(e_{i_{[g_k]}}) J_{[g_k]} |H|^{-1/2} \sum_{h\in H} \ol{\chi(h)} \xi_h \\ &=
 \sigma(e_{i_{[g_k]}}) |H|^{-1/2} \sum_{h\in H} \ol{\chi(h)} \xi_{g_k+h} \\ &=
 \alpha_0 |H|^{-1/2} \sum_{h\in H} \ol{\chi(h)} \xi_{g_k+\fg_1+h} \\ &=
 \alpha_0 |H|^{-1/2} \sum_{h\in H} \ol{\chi(h)} \xi_{g_l+h_{k1}+h} \\ &=
 \alpha_0 J_{[g_l]} \chi(h_{k1}) |H|^{-1/2}
   \sum_{h\in H} \ol{\chi(h+h_{k1})} \xi_{h+h_{k1}} \\ &=
 \alpha_0 \chi(h_{k1}) \zeta^\chi_{[g_l]} =
 \alpha_0 \chi(h_{k1}) \zeta^\chi_{[g_k+\fg_1]} .
\end{align*}
Similarly,
\[
 \sigma_\chi(f_{j_{[g_k]}}) \zeta_{[g_k]}^\chi
 = \beta_0 \chi(h_{k2}) \zeta^\chi_{[g_k+\fg_2]} .
\]

So the representations $\sigma_\chi$ all act on $\ltwo(G/H)$
with the same functions $i,j$, but with different constants.
Since $H$ is finite, so is $G$; and we may write $G = \bZ^2/K$
where $K$ is a subgroup of finite index.
Then $G/H = \bZ^2/HK =: \bZ^2/L$.
We wish to calculate the character $\psi_\chi$ on $L$
which distinguishes $\sigma_\chi$.

\begin{lem}
Let $\sigma$ be a group construction representation on
$\ltwo(\bZ^2/K)$ with symmetry group $L/K$,
and scalars determined by the character $\psi\in\hat{K}$.
Then the summands $\sigma_\chi$ are determined
by the set of characters $\psi_\chi\in\hat{L}$ satisfying $\psi_\chi|_K=\psi$.
The enumeration is given by elements of $\hat{H}=\widehat{L/K}$
which are related by $\psi_\chi(l) = \psi_0(l) \chi(l+K)$; and
this enumerates all possible extensions of $\psi$ from $K$ to $L$.
\end{lem}

\begin{pf}
For each $l \in L$, there is a unique word $w_l \in \Fth$ of degree $l$ so that
$\sigma_\chi(w_l) \dot \xi_0^\chi = \dot \xi_0^\chi$; and then
\[ \sigma_\chi(w_l) \xi_0^\chi = \psi_\chi(l) \xi_0^\chi .\]
When $k \in K$, one has
\[ \sigma(w_k) \xi_h = \psi(k) \xi_h \qforal h \in H.\]
Therefore it follows that
\[ \sigma_\chi(w_k) \xi_0^\chi = \psi(k) \xi_0^\chi \qforal h \in H.\]
So $\psi_\chi|_K=\psi$.

In general, we fix an extension $\psi_0$ of $\psi$ to $\bZ^2$
and use this to calculate $\psi_\chi$ as in Theorem~\ref{T:scalars}.
For $l \in L$, let $h_l := l+K \in L/K=H$.
\begin{align*}
 \psi_\chi(l) &= \frac1{|H|}
 \ip{ \sigma(w_l) \sum_{h_1\in H} \ol{\chi(h_1)} \xi_{h_1},
 \sum_{h_2\in H} \ol{\chi(h_2)} \xi_{h_2} } \\ &=
 \frac1{|H|} \sum_{h_1\in H} \sum_{h_2\in H}
 \ol{\chi(h_1)} \chi(h_2) \psi_0(l) \ip{\xi_{h_1+ h_l},\xi_{h_2}} \\ &=
 \chi(h_l) \psi_0(l).
\end{align*}
Hence we obtain the desired relationship between $\psi_0$ and $\psi_\chi$.

Note that the definition of the subspaces $\M_\chi$ depends on the
choice of $\phi$, but that the decomposition is unique.
All extensions of $\psi$ occur in this manner, so we have
enumerated all possibilities.
\bx\end{pf}

Lastly we explain why these summands are irreducible.
Clearly, if $H \ne \{0\}$, the representation is reducible
because we have exhibited a non-trivial collection of reducing subspaces.
The restriction to each of these subspaces yields a representation
of $G/H$ on $\ltwo(G/H)$.
By construction, the symmetry group of each of these representations
is $H/H = \{0\}$.
So irreducibility follows from the following lemma.

\begin{lem}\label{L:irreducible}
If a group construction representation of $\Fth$ on $\ltwo(G)$ has symmetry
group $\{0\}$, then it is irreducible.
\end{lem}

\begin{pf}
The complete lack of symmetry means that for any
element $g\in G\bsl\{0\}$, there is a word $w\in\Fth$ so that
$\sigma(w)\xi_0 \ne 0$ and $\sigma(w)\xi_g = 0$.
Indeed, if we write $G = \{ g_0=0,g_1,g_2,\dots\}$,
there are words $w_k$ so that
\[
 \sigma(w_k)\xi_0 \ne 0 \qand
 \sigma(w_k)\xi_{g_i} = 0 \FOR 1 \le i \le k .
\]
Thus $P_k = \sigma(w_k)^* \sigma(w_k)$ are projections
such that $P_k \xi_0=\xi_0$ for all $k\ge1$ and $P_k \xi_{g_i} = 0$
when $k \ge i$.  Therefore $\wotlim P_k = \xi_0\xi_0^*$.

For any $g_i\in G$, there are words $w_i,x_i \in \Fth$
so that $\sigma(w_i)^* \sigma(x_i) \xi_0 = \xi_{g_i}$.
Set $V_i := \sigma(w_i)^* \sigma(x_i)$.
Then $V_i \xi_0\xi_0^* V_j^* = \xi_{g_i} \xi_{g_j}^*$.
This is a complete set of matrix units for the compact operators
in the von Neumann algebra generated by $\sigma(\Fth)$.
Therefore $\sigma$ is irreducible.
\bx\end{pf}

A similar analysis works in the case of infinite groups.
Of course, there are no subspaces corresponding to the
representations $\sigma_\chi$; but the procedure is just a
measure theoretic version of the same.
Putting all of this together, we obtain the following decomposition.

\begin{thm} \label{T:decomposition}
Let $\sigma$ be a group construction representation of $\Fth$ on $G = \bZ^2/K$.
If the symmetry group $H$ is finite, then $\sigma$ decomposes as a
direct sum of irreducible atomic representations $\sigma_\chi$
on $\ltwo(G/H)$, one for each $\chi \in \hat{H}$.
If $H$ is infinite, then $\sigma$ decomposes as a direct integral
over the dual group $\hat{H}$ of irreducible atomic representations
$\sigma_\chi$ on $\ltwo(G/H)$.
The representations $\sigma_\chi$ all have the same graph,
and the scalars are given by all possible extensions of the
character $\psi$ on $K$ to $L=HK$.
\end{thm}

\medbreak\textbf{Further considerations.}
Because of the Dilation Theorem~\ref{defectfreediln}, we know that an
atomic $*$-dilation is uniquely determined by its restriction
to any cyclic coinvariant subspace.
In the cases examined in this section, one can often select a
smaller subspace which will suffice.

In the case of a finite group, the subspace has no proper subspace
which is cyclic and coinvariant.  But when the group is infinite, there
are many such subspaces---and none are minimal.
For example, in Example~\ref{3a_repn} we saw that the restriction
of a representation on $\bZ^2$ to $(-\bN_0)^2$ is such a subspace.
Indeed, if the subspace is spanned by standard basis vectors, then
whenever if contains $\xi_{s_0,t_0}$, it must contain all $\xi_{s,t}$
for $s \le s_0$ and $t \le t_0$.
Likewise, if $G = \bZ^2 / \bZ(p,q)$ with $pq \le 0$, then
$\spn\{ \xi_{[s,t]} : s \le s_0,\ t \le t_0 \}$ is a proper cyclic
coinvariant subspace.
In the case $pq>0$, it is easy to see that there is no proper
cyclic invariant subspace.  However we shall also see that
$\bZ(p,q)$ with $pq>0$ can never be the full symmetry group
of any representation of $\Fth$.

This discussion suggests that we need to put an equivalence relation
on these representations.  It is evident that two equivalent representations
must have unitarily equivalent $*$-dilations.

\begin{defn}
Two group construction representations $\sigma$ and $\sigma'$
of $\Fth$ on $G=\bZ^2/K$, with data $\{ \alpha_0, \beta_0, i_g,j_g: g \in G\}$
and $\{ \alpha'_0, \beta'_0, i'_g, j'_g : g \in G\}$ respectively, are
said to be \textit{equivalent}
if there is an integer $T$ and a character $\chi \in \hat{G}$ so that
\begin{alignat*}{2}
 \alpha'_0 &= \chi(\fg_1) \alpha_0,&\quad
 \beta'_0 &= \chi(\fg_2) \beta_0,\\
 i'_g &= i_g \AnD j'_g = j_g &\quad\FORAL g &= [s,t],\ s,t \le T .
\end{alignat*}
\end{defn}

Conversely, a group-like construction
which is defined only on a coinvariant
cyclic atomic subspace $\spn\{ \xi_g : g= [s,t],\ s,t \le T \}$
may be extended to the full group, usually in many ways.
To see this, take the minimal $*$-dilation $\pi$.
There are two cases, $G=\bZ^2$ and $G=\bZ^2/\bZ(p,q)$
with $pq\le 0$.
In the first case, take the basis vector $\xi_{T,T}$ and any
infinite word $\dots f_{j_k}e_{i_k} \dots f_{j_T}e_{i_T}$.
Define $\xi_{k,k} = \pi(f_{j_{k-1}}e_{i_{k-1}} \dots f_{j_T}e_{i_T}) \xi_{T,T}$.
The minimal coinvariant subspace $\M_k$ generated by $\xi_{k,k}$
can be identified with an atomic basis $\{\xi_g : g= [s,t],\ s,t \le k \}$.
These subspaces are nested, and their union yields a defect free
atomic representation on $\ltwo(G)$ as desired.
In the case of $G=\bZ^2/\bZ(p,q)$ with $pq\le 0$, we may suppose that
$q \ne 0$.  Then using an infinite word $\dots e_{i_k} \dots e_{i_T}$
will work in exactly the same way.
\smallbreak

The other issue to discuss here is how the decomposition of the
$*$-dilation corresponds to the decomposition of the restriction
to the coinvariant subspace.  This is a direct consequence of the
uniqueness of minimal dilations.
Indeed in the case of a representation which decomposes into
a direct sum of irreducible representations, the direct sum of the
$*$-dilations of the summands is clearly a minimal $*$-dilation.
Hence it is the unique $*$-dilation.
That is, the $*$-dilation of a direct sum is the direct sum of
the $*$-dilations of the summands.

\begin{thm}
Let $\sigma$ be a group construction representation of $\Fth$ on $\ltwo(G)$
with symmetry group $H$.  Then the minimal $*$-dilation
decomposes as a direct integral of the $*$-dilations of the
irreducible integrands in the direct integral decomposition of $\sigma$.
\end{thm}

\begin{pf}
Let $\pi$ denote the minimal $*$-dilation of $\sigma$ on $\K$.
We first show that there is a spectral measure on $\K$ over
measurable subsets of $\hat{H}$ which is absolutely continuous
with respect to Haar measure and extends the spectral measure
on $\ltwo(G)$.

Write $\sigma$ as a direct integral over $\hat{H}$.
For each measurable subset $A \subset \hat{H}$, let
$E(A)$ be the spectral projection onto $\spn\{ \M_\chi : \chi \in A \}$.
The restriction $\sigma_A$ of $\sigma$ to this subspace
has a minimal $*$-dilation $\pi_A$.
By uniqueness, $\pi \simeq \pi_A \oplus \pi_{A^c}$.
This decomposition splits $\K \simeq F(A)\K \oplus F(A^c)\K$.
It is routine to check that $F$ is countably additive, that
$E(A) = P_\H F(A) = F(A) P_\H$, and that $F$ is absolutely continuous.

The rest follows from standard arguments.
Since we do not actually need any explicit formulae for the
direct integral decomposition for any of our analyses, we will not
subject the reader to the technicalities.
\bx\end{pf}

\medbreak\textbf{The Main Theorem.}
The central result of this paper is the following:

\begin{thm} \label{T:main}
Every atomic $*$-representation of $\Fth$ with connected graph is the
minimal $*$-dilation of a group construction representation.
It is irreducible if and only if its symmetry group is trivial.
In general, it decomposes as a direct sum or direct integral
of irreducible group construction representations.
\end{thm}

The proof evolves from a case by case analysis of the various types.

\section{Finitely Correlated Atomic Representations}
\label{S:atomicDC}

In this section, we restrict our attention to atomic
$*$-represent\-ations of $\Fth$ which are finitely correlated.
Such representations are particularly tractable.
As in the case of the free semigroup \cite{DKS}, the whole
class of finitely correlated $*$-representations
may turn out to be classifiable. This general problem
is not considered here.

We assume that the graph is connected.
By the discussion in the Section~\ref{S:atomic}, we obtain a coinvariant cyclic
subspace spanned by the ring by ring portion of the graph.
We wish to show that it always arises from a group construction.
To this end, we need a criterion for when this construction is possible.

\begin{lem}\label{L:commute}
Let $\rho$ be an atomic representation of $\Fth$.
If $\dot\xi$ is a vertex of $\G_\rho$ and $u,v$ are words such that
$ \rho(e_u) \dot\xi = \dot\xi = \rho(f_v) \dot\xi ,$
then $e_uf_v = f_ve_u$.
\end{lem}

\begin{pf}
There are words $u'$ and $v'$ so that $e_uf_v=f_{v'}e_{u'}$.
So $\rho(f_{v'}e_{u'}) \dot\xi = \dot\xi = \rho(f_v) \dot\xi$.
Since $|v'|=|v|$,  $\rho(f_{v'})$ would have range orthogonal
to $\rho(f_v)$ unless $v'=v$.
Therefore they must be equal.
Similarly $u'=u$.
\bx\end{pf}

\begin{lem}\label{L:ringbyring}
Suppose that $e_{u_0}f_{v_0} = f_{v_0} e_{u_0}$ where $|u_0|=k$ and $|v_0|=l$.
Then there is an atomic defect free representation
$\sigma$ on $\C_k \times \C_l$
with $\sigma(e_{u_0})\xi_0 = \xi_0 = \sigma(f_{v_0}) \xi_0$.

Arbitrary constants $\alpha,\beta\in\bT$ yield a
representation $\sigma_{\alpha,\beta}$ given by
$\sigma_{\alpha,\beta}(e_i) = \alpha \sigma(e_i)$ and
$\sigma_{\alpha,\beta}(f_j) = \beta \sigma(f_j)$.
Then $\sigma_{\alpha,\beta} \simeq \sigma_{\alpha',\beta'}$
if and only if $\alpha^k = \alpha^{\prime k}$
and $\beta^l = \beta^{\prime l}$.
\end{lem}

\begin{pf}
Write $u_0 = i_{k-1,0} \dots i_{0,0}$ and ${v_0}= j_{0,l-1} \dots j_{0,0}$.
The commutation relations show that there are unique words $u_t$ for
$0 \le t \le l$, so that
\[ f_{v_0} e_{u_0} = f_{j_{0,l-1}} \dots f_{j_{0,t}} e_{u_t} f_{j_{0,t-1}}
\dots f_{j_{0,0}} .\] Write $u_t = i_{k-1,t} \dots i_{s,t} \dots
i_{0,t}$ for $0 \le t \le l$ and note that $u_l = u_0$. Similarly,
there are unique words $v_s$ so that
\[ f_{v_0} e_{u_0} = e_{i_{k-1,0}} \dots e_{i_{s,0}} f_{v_s} e_{i_{s-1,0}}
\dots e_{i_{0,0}} .\] Write $v_s =  j_{s,l-1} \dots j_{s,0}$ for $0
\le s \le k$; and again one has $v_k = v_0$.

It follows from unique factorization that
\begin{align*}
 f_{v_0} e_{u_0} &=
 e_{i_{k\!-\!1,0}} \dots e_{i_{s\!+\!1,0}} f_{j_{s\!+\!1,l\!-\!1}} \dots f_{j_{s\!+\!1,t}}
e_{i_{s,t}}
  f_{j_{s,t\!-\!1}} \dots f_{j_{s,0}} e_{i_{s\!-\!1,0}} \dots e_{i_{0,0}} \\
 &=
 e_{i_{k\!-\!1,0}} \dots e_{i_{s\!+\!1,0}} f_{j_{s\!+\!1,l\!-\!1}} \dots f_{j_{s\!+\!1,t\!+\!1}}
e_{i_{s,t\!+\!1}}
  f_{j_{s,t}} \dots f_{j_{s,0}} e_{i_{s\!-\!1,0}} \dots e_{i_{0,0}}.
\end{align*}
Now cancellation shows that
\[
 f_{j_{s+1,t}} e_{i_{s,t}} = e_{i_{s,t+1}} f_{j_{s,t}}
 \qforal s \in \C_k \AnD t \in \C_l .
\]

These are the relations needed to allow the construction of
a homomorphism on $\ltwo(\C_k\times\C_l)$ as in Section~\ref{S:group}.
Namely,
\[
 \sigma_{\alpha,\beta}(e_{i_{s,t}}) \xi_{s,t} = \alpha \xi_{s+1,t} \qand
 \sigma_{\alpha,\beta}(f_{j_{s,t}}) \xi_{s,t} = \beta \xi_{s,t+1} .
\]
One calculates that $\sigma_{\alpha,\beta}(e_{u_0}) \xi_{0,0} = \alpha^k \xi_{0,0}$ and
$\sigma_{\alpha,\beta}(f_{v_0}) \xi_{0,0} = \beta^l \xi_{0,0}$.

Clearly $\sigma_{\alpha,\beta}$ is also a defect free
atomic representation for any $\alpha,\beta \in \bT$.
Indeed, scalars can be assigned to each edge arbitrarily;
but by Theorem~\ref{T:scalars}, there is no loss in making the
constants all the same.
The characters of $\C_k \times \C_l$ have the form
$\chi(s,t) = \omega_1^s\omega_2^t$ where
$\omega_1^k = 1 = \omega_2^l$.
Thus Theorem~\ref{T:scalars} also shows that
$\sigma_{\alpha,\beta} \simeq \sigma_{\alpha',\beta'}$
if and only if $\alpha' = \chi(\fg_1) \alpha = \omega_1\alpha$
and $\beta' = \chi(\fg_2) \beta = \omega_2 \beta$; and this
is equivalent to  $\alpha^{\prime k} = \alpha^k$
and $\beta^{\prime l} = \beta^l$.
\bx\end{pf}

\begin{cor}
Let $H$ be the symmetry subgroup for the representation $\sigma_{\alpha,\beta}$
constructed in Lemma~$\ref{L:ringbyring}$.
For any subgroup $K \le H$, there is a group construction representation on
$\ltwo(G/K)$ such that
\[
 \sigma(e_{u_0})\xi_0 = \alpha^k \xi_0
 \qand \sigma(f_{v_0})\xi_0 = \beta^l \xi_0 .
\]
This representation decomposes as a direct sum of irreducible
atomic defect free representations on $\ltwo(G/H)$ indexed
by $\widehat{H/K}$.
\end{cor}

\begin{pf}
It is evident that the induced representation on $\ltwo(G/K)$
fits the conditions of Section~\ref{S:group}.
The symmetry group is clearly $H/K$.
Since $H/K$ is finite, the decomposition Theorem~\ref{T:decomposition}
yields a finite direct sum of irreducible representations
on $\ltwo(G/H)$.
\bx\end{pf}

We are now ready to establish the following result.

\begin{thm}\label{t.f.cor.d.c.}
Any defect free atomic representation
of $\Fth$ on a finite dimensional space with connected graph is
isometrically isomorphic to a dilation of a group construction
representation for a finite group $G$.
\end{thm}

\begin{pf}
Let $\sigma $ be a given finitely correlated defect free
representation on a finite dimensional Hilbert space $\H$,
and with a connected graph.
Pulling back along the blue (respectively red) edges eventually
reaches a periodic state in the ring by ring portion of the
representation.
Fix a standard basis vector $\xi_0$ in the ring by ring.
One can find the unique minimal words $u_0$ and $v_0$
so that $e_{u_0} \dot\xi_0 = \dot\xi_0$ and
$f_{v_0} \dot\xi_0 = \dot\xi_0$.
Then $e_{u_0} f_{v_0} = f_{v_0} e_{u_0}$ by Lemma~\ref{L:commute}.

For each $(s,t) \in (\bN_0)^2$, there is a unique word $w_{s,t} \in \Fth$
of degree $(s,t)$ such that $\sigma(w_{s,t})\dot\xi_0 =: \dot\xi_{s,t} \ne 0$.
In particular, $w_{k,0} = e_{u_0}$ and $w_{0,l} = f_{v_0}$.
From the uniqueness of this path,
$w_{s+ak,t+bl} = w_{s,t}e_{u_0}^af_{v_0}^b$.
So $\dot\xi_{s+ak,t+bl} = \dot \xi_{s,t}$;
thus the set of vectors $\xi_{s,t}$ is periodic, and so may be indexed
by an element $g$ of $\C_k\times\C_l$, say $\xi_g$.
Let $K$ denote the set
\[ K = \{g \in \C_k\times\C_l : \dot\xi_g = \dot\xi_0 \} .\]
Clearly $K$ is closed under addition, and hence is a
subgroup of $\C_k\times\C_l$.
The cosets each determine a distinct basis vector.
So the subspace $\H_0$ spanned by all the basis vectors $\xi_{s,t}$
in the ring by ring portion of the graph is naturally identified
with $\ltwo(G)$, where $G = \C_k\times\C_l/K$.
Observe that $\H_0$ is a cyclic coinvariant subspace for $\sigma$.

Thus there is a setup exactly as in section~\ref{S:group}
using the group $G$, with the functions $i$, $j$, $\alpha$
and $\beta$ determined by the action of $\sigma$ on the
basis $\{\xi_g : g \in G \}$ for $\H$.
The compression $\sigma'$ of $\sigma$ to $\H_0$ is therefore
unitarily equivalent to this group construction.
The representation $\sigma$ has a unique minimal $*$-dilation $\pi$,
and this evidently coincides with the minimal $*$-dilation of $\sigma'$.
Hence there is a subspace of the $*$-dilation of $\sigma'$
which corresponds to $\H$, exhibiting $\sigma$ as a dilation of $\sigma'$.
\bx\end{pf}

By Theorem \ref{t.f.cor.d.c.} one can now easily obtain the following.
It is sufficient to decompose the graph into connected components.

\begin{cor}\label{f.cor.d.c.}
Any finitely correlated atomic $*$-rep\-re\-sen\-ta\-tion
of $\Fth$ is unitarily equivalent to the direct sum of
irreducible atomic $*$-repre\-sent\-ations which dilate
group construction representations on finite abelian groups.
\end{cor}

\begin{exmp}\label{cycle}
Let us take another look at Example~\ref{E:cycle reps}.
Let $\theta$ be a permutation in $S_{mn}$, and fix a cycle of $\theta$:
\[
 \big( (i_0,j_0), (i_1,j_1), \dots, (i_{k-1},j_{k-1}) \big).
\]
Let $G = \C_k$ with $\fg_1 = 1$ and $\fg_2 = -1$.
Set $i(g) = i_g$ and $j(g) = j_{g-1}$ for $g \in G$.
Use this to define a representation $\rho_{\alpha,\beta}$.
Then at each vertex $\xi_g$,
\begin{align*}
 \rho_{\alpha,\beta}(f_{j_g} e_{i_g}) \xi_g
 &= \alpha \rho_{\alpha,\beta}(f_{j_g}) \xi_{g+1}
 = \alpha\beta\xi_g\\
 &= \beta \rho_{\alpha,\beta}(e_{i_{g-1}}) \xi_{g-1}
 = \rho_{\alpha,\beta}(e_{i_{g-1}} f_{j_{g-1}}) \xi_g .
\end{align*}
Thus the commutation relations show that this is a representation.
It is not difficult to see that there are no symmetries.
So it is irreducible.

Let $u_0=i_{k-1}\dots i_1i_0$ and $v_0=j_0j_1\dots j_{k-1}$.
These elements satisfy the identities
$\rho_{\alpha,\beta}(e_{u_0})\xi_0 = \alpha^k \xi_0$
and $\rho_{\alpha,\beta}(f_{v_0})\xi_0 = \beta^k \xi_0$.
It is easy to check that $e_{u_0}f_{v_0}=f_{v_0}e_{u_0}$.
Therefore we may consider the atomic representation $\sigma$
on $\C_k \times \C_k$ given by Lemma~\ref{L:ringbyring}.

It is easy to obtain, for $0 \le s < k$, that
\begin{align*}
 u_s &= i_{k-1-s} \dots i_1 i_0 i_{k-1}\dots i_{k-s}
 \qand
 v_s = j_s j_{s+1}\dots j_{k-1}j_0\dots j_{s-1} .
\end{align*}
It follows that the subgroup of symmetries includes
\[
H = \{ (0,0),(1,1),..., (k-1,k-1) \} .
\]
However, as this comes from a cycle of $\theta$, a little thought
shows that there are no other symmetries.  Compute
\[
 \C_k \times \C_k/H \simeq \C_k .
\]

The characters of $\C_k$ are given by $\chi(1)=\omega$
where $\omega^k = 1$.
Therefore by Theorem~\ref{T:decomposition}, we can decompose
$\sigma \simeq  \bigoplus_{\omega^k=1} \rho_{\omega,\ol{\omega}}.$
\end{exmp}

\medbreak\textbf{Long commuting words.}
We now show that there are many infinitely many finitely correlated representations
for any $\Fth$ by exhibiting arbitrarily long primitive commuting words.

\begin{prop}\label{LCW}
For any given $\Fth$, there are commuting pairs
$e_u$ and $f_v$ which determine irreducible atomic
representations of arbitrarily large dimension.
\end{prop}

\begin{pf}
Consider the two sets of words $\U:=\{e_u:|u|=N!\}$
and $\V:=\{f_v:|v|=N!\}$. Since $|\U|=m^{N!}$ and $|\V|=n^{N!}$, for
simplified notation, we write $\U=\{u_1,...,u_{m^{N!}}\}$ and
$\V=\{v_1,...,v_{n^{N!}}\}$.

Given $u_i\in\U$ and $v_j\in\V$, the relation $\theta$ uniquely
determines $u_{i'}\in\U$ and $v_{j'}\in\V$ such that
$e_{u_i}f_{v_j}=f_{v_{j'}}e_{u_{i'}}$.
We obtain a permutation $\theta'$ of $m^{N!}\times n^{N!}$
so that $\theta'(u_i,v_j) = (u_{i'},v_{j'})$.

If $e_uf_v=f_ve_u$ for all $u\in \U$ and $v\in\V$,
then pick any primitive words $u\in\U$ and $v\in\V$.
Then $e_uf_v= f_ve_u$ leads to an atomic representation on
a quotient $G$ of $\C_{N!} \times \C_{N!}$ by the
symmetry subgroup $H$.
Since $u$ has no symmetries, $(s,0)\notin H$ for any $1\le s<|u|$.
So $|G|\ge N!$.

Otherwise, $\theta'$ has a cycle, say $\C$, of length $t\ge 2$.
Without loss of generality, we can assume that
$\C=\big( (u_0,v_0),...,(u_{t-1},v_{t-1}) \big)$.
As in Example~\ref{cycle}, we obtain a representation $\rho$
on $\ltwo(\C_{N!t})$.
Therefore $N!$ does not belong to the symmetry subgroup $H$.
But $H = \ip{p}$ where $p=0$ or $p$ divides $N!t$.
The $t$ pairs $(u_i,v_i)$ are all distinct, and so at least one
of the words $U$ or $V$ is not $N!$-periodic.
So $p$ does not divide $N!$.
Therefore $p\ge N+1$ implying $|\C_{N!t}/H| \ge N+1$.
Consequently there is an irreducible representation on a space
of finite dimension greater than $N$.
\bx\end{pf}

\begin{exmp}\label{E:commuting favourite}
Consider the forward 3-cycle algebra of
Example~\ref{E:favourite} given by the permutation
$\big( (1,1) , (1,2), (2,1) \big)$ in $S_{2\times 2}$.
The relations have the succinct form $f_je_i = e_{i+j} f_i$
where addition is modulo 2. Observe that
\begin{gather*}
 f_{i_k} e_{i_1,i_2, \dots, i_k} =
 e_{i_k\!+i_1,i_1\!+i_2,\dots,i_{k\!-\!1}\!+i_k} f_{i_k}  \\
 f_{i_{k\!-\!1}\!+i_k}
 e_{i_k\!+i_1,i_1\!+i_2,\dots,i_{k\!-\!1}\!+i_k} =
 e_{i_{k\!-\!1}\!+2i_k\!+i_1,i_k
 \!+2i_1\!+i_2,\dots,i_{k\!-\!2}\!+2i_{k\!-\!1}\!+i_k}
 f_{i_{k\!-\!1}\!+i_k} \\
 \quad \vdots\\
 f_{\sum\! \binom{n\!-\!1}p i_{k\!-\!p}} e_{\sum\! \binom{n\!-\!1}p
 i_{1\!-\!p},\dots,\sum\! \binom{n\!-\!1}p i_{k\!-\!p}} = e_{\sum\!
 \binom n p i_{1\!-\!p},\dots,\sum\! \binom n p i_{k\!-\!p}}
 f_{\sum\! \binom{n\!-\!1}p i_{k\!-\!p}}
\end{gather*}
The binomial sums are over $p\ge 0$. If we take $k=l=2^n-1$, then
$\binom{2^n-1}p \equiv 1$ for all $p$. So in order that the last $e$
term equal the original, it suffices that $\sum_{s=1}^l i_s \equiv 0$.
The $f$ word is then uniquely determined so that $f_ve_u = e_uf_v$,.

For example,
\[ f_{1222212} e_{1121212} = e_{1121212} f_{1222212} .\]
It follows that there are primitive commuting pairs of arbitrarily
great length.
\end{exmp}

\section{The Ring by Tail case}
\label{S:ringbytail}

Now consider Case 2a, the ring by tail type in which the blue components are ring
representations and the red components are infinite tail representations.

To simplify the presentation, we introduce some notation.
For a word $u = i_{k-1} \dots i_{0}$ of length $k\ge1$,
define
\[
 u(s,0] := i_{s-1}\dots i_{0} \qand u^{(s)} := i_{s-1}\dots i_{0} i_{k-1} \dots i_{s}
 \qfor 0\le s < k.
\]
and observe that $u^{(k)} = u$.
Also if $v = v_{-1}v_{-2}v_{-3}\dots$ is an infinite word, for $0 \le t < t'$ write
\[
 v(-t,-t'] = v_{-t-1}\dots v_{-t'} \qand
 v(-t,-\infty) =  v_{-t-1}v_{-t-2}v_{-t-3}\dots.
\]

Start with a basis vector, $\xi_{0,0}$ say,
in one of the blue rings of minimal size $k$.
One can find a unique word $u_0 = i_{k-1,0} \dots i_{0,0}$
so that $\rho(e_{u_0}) \dot\xi_{0,0} = \dot\xi_{0,0}$.
Choose a scalar $\alpha$ so that
\[ \rho(e_{u_0})  \xi_{0,0} =\alpha^k \xi_{0,0} .\]
For $1 \le s \le k-1$, renormalize the basis vectors so that
$\xi_{s,0} = \ol{\alpha}^s \rho(e_{u_0(s,0]}) \xi_{0,0}$.
Then
\[ \rho(e_{i_{s,0}})  \xi_{s,0} =\alpha \xi_{s+1,0} \qfor s \in \C_k .\]

Pull back along the red edges from $\xi_{0,0}$ to obtain
an infinite red tail $v_0 :=j_{0,-1}j_{0,-2}\dots$ and
basis vectors $\dot\xi_{0,t}$ for $t < 0$ satisfying
\[ \rho(f_{j_{0,t}}) \dot \xi_{0,t} = \dot\xi_{0,t+1} \qfor t<0 .\]
Normalize the basis vectors so that
$ \rho(f_{j_{0,-t}}) \xi_{0,-t} = \xi_{0,1-t}$ for $t\ge1$.

There is a unique word $u_t = i_{k-1,t} \dots i_{0,t}$
of length $k$ for each $t<0$ so that
$e_{u_0} f_{v_0(0,t]} = f_{v_0(0,t]}   e_{u_{t}}$.
It follows that $\rho(e_{u_{t}}) \dot\xi_{0,t} = \dot \xi_{0,t}$;
and hence one can deduce that
$\rho(e_{u_{t}}) \xi_{0,t} = \alpha^k \xi_{0,t}$ for $t<0$.
Define $\xi_{s,t} = \ol{\alpha}^s \rho(u_{t}(s,0]) \xi_{0,t}$.
Then
\[ \rho(e_{i_{s,t}})  \xi_{s,t} =\alpha \xi_{s+1,t} \qfor s \in \C_k \AnD t \le 0.\]

Similarly one can pull back along the red edges from $\xi_{s,0}$
to obtain an infinite word $v_s :=j_{s,-1}j_{s,-2}\dots$ for $1 \le s \le k-1$.
Again using the commutation relations, one obtains that
\[ \rho(f_{j_{s,t}}) \xi_{s,t} =\xi_{s,t+1} \qfor t<0 \AnD s \in \C_k.\]
It is also easy to verify that the cycles $u_{t}^{(s)}$ satisfy
\[ \rho(e_{u_{t}^{(s)}}) \xi_{s,t} = \alpha^k \xi_{s,t}  \qfor t\le 0 \AnD s \in \C_k.\]

It is evident that the vectors $\{ \xi_{s,t} : 0 \le s < k \AnD t \le0 \}$
span a coinvariant subspace.
By the connectedness of the graph, it is also a cyclic subspace.
Thus this subspace determines the representation by the uniqueness of the
isometric dilation in Theorem~\ref{defectfreediln}.

There are only $m^k$ words of length $k$ in $m$ letters.
So by the pigeonhole principle, there is some word $u$ which is
repeated infinitely often in the sequence $\{u_{t} : t \le0 \}$.
Without loss of generality, we may assume that there is a sequence
$t_0=0>t_1>t_2>\dots$ such that $u_{t_k} = u_0$ for $k \ge1$.
Then
\[ e_{u_0} f_{v_0(0,t_k]} = f_{v_0(0,t_k]} e_{u_0} \qfor k\ge1 .\]

Conversely, given a word $u_0$, $\alpha \in \bT$,
an infinite tail $v_0 = j_{0,-1}j_{0,-2}\dots$ and
a sequence $0<t_1<t_2<\dots$ such that
$e_{u_0} f_{v_0(0,-t_k]} = f_{v_0(0,-t_k]} e_{u_0}$ for $k\ge1$,
one can build a representation of type Case 2a.
We will do this by our group construction by extending
the definition of this subspace indexed by $\C_k \times -\bN_0$
to the group $G=\C_k \times \bZ$.
This may be accomplished in many ways, but a simple way is to
make it periodic on $\C_k \times \bN_0$ by repeating the
segment on $\C_k \times [0,t_1)$ using the fact that
$u_{t_1}=u_0$ and $e_{u_0} f_{v_0(0,t_1]} = f_{v_0(0,t_1]} e_{u_0}$.
The finitely correlated representation of Lemma~\ref{L:ringbyring}
can be unfolded to obtain a representation on $\C_k\times \bZ$.
Just cut off the left half of the cylinder and glue it to the
one we have.
Any choice of extension is tail equivalent to any other.
\smallbreak

Note that this analysis is necessary, and it is not the
case that an arbitrary tail $e_{u_0}$ and infinite word
$f_{v_0}$ determines a representation.  For example,
this may force a red edge and blue edge into the same
vertex which is not possible from the commutation relations.

We have obtained:

\begin{thm}\label{T:ringbytail}
A ring by tail representation is determined by a word $u_0$
of length $k$, a scalar $\alpha \in \bT$,
an infinite word $v_0 :=j_{0,-1}j_{0,-2}\dots$ and
a sequence $0>t_1>t_2>\dots$ such that
$e_{u_0} f_{v_0(0,t_k]} = f_{v_0(0,t_k]} e_{u_0}$ for $k\ge1$.
It corresponds to a group construction on the group $\C_k \times \bZ$.
\end{thm}

Since a ring by tail representation is given by the group
construction on $G=\C_k \times \bZ$, we can use the
decomposition results of Section~\ref{S:group}.

\begin{thm}\label{T:ringbytailreduction}
Let $\pi$ be a (connected) ring by tail representation
with symmetry subgroup $H \le \C_k \times \bZ$.
If $H  \le \C_k \times \{0\}$, then $\pi$ decomposes as
a finite direct sum of irreducible ring by tail representations.
Otherwise, it decomposes as a direct integral of irreducible
ring by ring atomic representations.
\end{thm}

\begin{pf}
In the first case, $H = \ip{(d,0)}$ is a subgroup of $\C_k \times \{0\}$.
Since $H$ is finite, Theorem~\ref{T:decomposition} splits
$\pi$ as a direct sum of finitely many irreducible representations
for the quotient group $G/H \simeq \C_{k/d}\times\bZ$.
So these are irreducible ring by tail representations.

In the second case, $H$ contains an element $(a,b)$ with $b \ne 0$.
Therefore $H$ contains the element $(0,kb)$.
Thus $G=\C_k\times\bZ /H$ is a quotient of $\C_k \times \C_{kb}$;
and in particular, $G$ is finite.
Since $H$ is infinite, Theorem~\ref{T:decomposition} yields
a direct integral of irreducible representations on $\ltwo(G)$.
Evidently, these are ring by ring type.
\bx\end{pf}

\smallskip
Case 2b is handled by exchanging the role of the red and blue edges in Case 2a.

\section{The Tail by Tail case}
\label{S:tailbytail}

In Case 3, both the red and the blue components are infinite tail representations.
The most basic is Case 3a, in which the red tail never intersects the original
blue component again.

\medbreak\textbf{Case 3a.}
As we saw in Example~\ref{3a_repn}, the 3a case is an inductive limit
of copies of the left regular representation.
So in a certain sense, they are the easiest.

Given a representation of type 3a,
start with any basis vector $\xi_{0,0}$.
Pull back on both blue and red edges to determine integers $i_{s,t}$
and $j_{s,t}$ and basis vectors  $\dot\xi_{s,t}$ for $s,t \le 0$ so that
\begin{align*}
 \pi( e_{i_{s,t}}) \dot\xi_{s-1,t} &= \dot\xi_{s,t} \qand \\
 \pi( f_{j_{s,t}}) \dot\xi_{s,t-1} &= \dot\xi_{s,t}  \qfor s, t \le 0 .
\end{align*}
The assumption of Case 3a ensures that the $\dot\xi_{s,t}$ are all distinct.
Thus we have found a cyclic coinvariant subspace $\M = \spn\{\xi_{s,t} : s,t \le 0\}$
on which we have a representation on $(-\bN_0)^2$.

The representation is determined by the tail
\[ \tau = e_{0,0}f_{-1,0}e_{-1,-1}f_{-2,-1} \dots \]
since the other $i_{s,t}$ and $j_{s,t}$ for $s,t \le 0$
are determined by the commutation relations.
One could use Theorem~\ref{T:scalars} to make the scalars
all equal to 1; but in fact that is not necessary.
Instead one selects the appropriate unit vector
$\xi_{s,t}$ in $\dot\xi_{s,t}$ recursively so that
\begin{align*}
 \pi(e_{i_{s,t}}) \xi_{s-1,t} &= \xi_{s,t} \qand \\
 \pi( f_{j_{s,t}}) \xi_{s,t-1} &= \xi_{s,t}  \qfor s, t \le 0 .
\end{align*}

Example~\ref{3a_repn} explains how to dilate this to a $*$-rep\-re\-sen\-ta\-tion
which is an inductive limit of copies of the left regular representation.
By Theorem~\ref{defectfreediln}, this is the unique minimal $*$-dilation;
and hence it is $\pi$.

While the infinite tail is sufficient to describe the representation,
as is done in Example~\ref{3a_repn},
it is not unique, and the equivalence relation on tails is
not at all transparent.
A much more useful collection of data associated to $\tau$
is the set
\[
 \Sigma(\pi,\xi_{0,0}) = \Sigma(\tau) = \{ (i_{s,t}, j_{s,t}) : s, t \le0 \} .
\]
In Definition~\ref{D:shift tail}, we put an equivalence relation
of shift tail equivalence on these sets.

\begin{defn}
For each inductive $*$-representation $\pi$, define
$\Sigma(\pi)$ to be the equivalence class of
$\Sigma(\pi,\xi)$ modulo shift tail equivalence
for any standard basis vector $\xi \in \H_\pi$.
\end{defn}

That this definition makes sense is part of the following result.

\begin{thm}
If $\pi$ is an inductive (type 3a) atomic $*$-rep\-re\-sen\-ta\-tion,
then $\Sigma(\pi)$ is independent of the choice of initial vector.
Two inductive $*$-rep\-re\-sen\-ta\-tions $\pi_1$ and $\pi_2$ are
unitarily equivalent if and only if $\Sigma(\pi_1) = \Sigma(\pi_2)$.
\end{thm}

\begin{pf}
Start with two standard basis vectors, $\xi_{0,0}$ and $\zeta_{0,0}$.
The connectedness of the graph means that there is a path
from $\xi_{0,0}$ to $\zeta_{0,0}$.
By Lemma~\ref{L:push pull}, there is a path from $\xi_{0,0}$ to $\zeta_{0,0}$
of the form $uv^*$.
Let $d(v)=(s_1,t_1)$ and $d(u)=(s_2,t_2)$.
Then $\dot\xi_{-s_1,-t_1} = \dot\zeta_{-s_2,-t_2}$.
Therefore the data agrees on all basis vectors obtained
by pulling back from that common vector.
So the two data sets are $(s_1-s_2,t_1-t_2)$-shift tail equivalent.

Clearly, if two inductive $*$-representations $\sigma_1$ and $\sigma_2$
have shift tail equivalent data $\Sigma(\pi_i, \xi_i)$,
then they are unitarily equivalent.
Consider the converse.
Suppose that $\pi_1$ and $\pi_2$ are unitarily equivalent,
say via a unitary $U$ in $\B(\H_{\pi_1},\H_{\pi_2})$.
Fix the basis vector $\xi_{0,0}$ for $\pi_1$, and corresponding basis
$\xi_{s,t}$ for $s,t \le 0$ and data $\Sigma(\pi_1,\xi_{0,0})$.
 From the unitary equivalence, $\xi_{0,0}$ is identified with
a vector $\eta_{0,0} = U\xi \in \H_{\pi_2}$.
Write $\eta$ in the standard basis for $\pi_2$,
say $\eta = \sum a_i \zeta_i$; and choose a standard basis
vector $\zeta= \zeta_{i_0}$ for which  $a_{i_0} \ne 0$.
For any $s,t \le 0$, there is a word $w_{s,t}\in\Fth$
with $d(w_{s,t})=(|s|,|t|)$ so that $\pi_1(w_{s,t})\xi_{s,t} = \xi_{0,0}$.
Therefore $\eta_{0,0} = \pi_2(w_{s,t}) U\xi_{s,t}$ is in the range
of the partial isometry $\pi_2(w_{s,t})$.
Since $\pi_2$ is atomic, the range of $\pi_2(w_{s,t})$ is spanned
by standard basis vectors.
Consequently $\zeta$ is in the range of $\pi_2(w_{s,t})$ for
every $s,t \le 0$.
Restating this another way, it says that
$\Sigma(\pi_2,\zeta) = \Sigma(\pi_1,\xi_{0,0})$.
Hence $\Sigma(\pi_1) = \Sigma(\pi_2)$.
\bx\end{pf}

\begin{exmp}
It is tempting to think that an inductive $*$-rep\-re\-sen\-ta\-tion can be
determined by the blue and red components containing
a standard basis vector.
That is, start at a basis vector $\xi_{0,0}$.
Pull back along the blue edges to obtain the infinite word
$u = i_{-1}i_{-2}\dots$ determining the blue infinite tail
component  containing $\xi_{0,0}$.
Likewise pull back along the red edges to get a red infinite
tail $v = j_{-1}j_{-2}\dots$ determining the red component
containing $\xi_{0,0}$.
\textit{The infinite words $u$ and $v$ may not determine
the representation uniquely!}

Consider the permutation $\big( (1,1) , (1,3) \big) \big( (1,2) , (2,1) \big)$
in $S_{3 \times 3}$ with two 2-cycles;
and the two infinite words $u = e_1 e_3  e_3 e_3 \dots$
and $v = f_1 f_3 f_3 f_3 \dots$.

Two inequivalent $*$-rep\-re\-sen\-ta\-tions $\pi_1$ and $\pi_2$ will
be constructed to produce a coinvariant subspace spanned
by vectors $\xi_{-s,-t}$ for $s,t \ge 0$ with
\begin{alignat*}{2}
\pi_i(e_1) \xi_{-1,0} &= \xi_{0,0} &\qquad \pi_i(e_3) \xi_{-s-1,0} &= \xi_{-s,0} \qfor s \ge1\\
\pi_i(f_1)  \xi_{0,-1} &= \xi_{0,0} &\qquad \pi_i(f_3)  \xi_{0,-t-1} &= \xi_{-0,-t} \qfor t \ge1.
\end{alignat*}
So they have the same infinite blue and red components.
We define in addition
\begin{alignat*}{2}
\pi_1(e_3) \xi_{-s-1,-t} &= \xi_{-s,-t} &\qfor & s \ge 1 \AnD t \ge 0\\
\pi_1(e_2) \xi_{-1,-t} &= \xi_{0,-t} &\qfor &  t \ge 1\\
\pi_1(f_3) \xi_{-s,-t-1} &= \xi_{-s,-t} &\qfor & s \ge 0 \AnD t \ge 1\\
\pi_1(f_2) \xi_{-s,-1} &= \xi_{-s,0} &\qfor &  s \ge 1
\end{alignat*}
and
\begin{alignat*}{2}
\pi_2(e_3) \xi_{-s-1,-t} &= \xi_{-s,-t} &\qfor & s \ge 1 \AnD t \ge 0\\
\pi_2(e_1) \xi_{-1,-t} &= \xi_{0,-t} &\qfor &  t \ge 1\\
\pi_2(f_1) \xi_{-s,-t-1} &= \xi_{-s,-t} &\qfor & s \ge 0 \AnD t \ge 1\\
\pi_2(f_3) \xi_{-s,-1} &= \xi_{-s,0} &\qfor &  s \ge 1
\end{alignat*}
These two $*$-rep\-re\-sen\-ta\-tions are evidently not shift tail equivalent.
Thus they are not unitarily equivalent.

This means that these representations are not some form
of twisted product of an infinite tail representation for $\A_m$
with an infinite tail representation for $\A_n$.
\end{exmp}

\medbreak\textbf{Case 3b.} Now suppose that we have a $*$-rep\-re\-sen\-ta\-tion
$\pi$ of type 3b.
Pulling back along the red edges yields a periodic sequence of blue
components of infinite tail type, say $H_0, \dots,H_{l-1}$,
so that every red edge into each $H_i$ comes from $H_{i-1\!\! \pmod l}$.

Start with a basis vertex $\dot\xi$.
Let $v$ be the unique word of length $l$ so that $f_v$ maps onto $\dot\xi$.
Then there is a vertex $\dot\zeta \in H_0$ so that
$\pi(f_v) \dot\zeta = \dot\xi$.
Our first goal is to explain why these two vertices are comparable in $H_0$
if $\dot\xi$ is sufficiently far up the tail.

\begin{lem}\label{L:comparable}
Suppose that $\pi$ is a type {\em 3b} atomic $*$-rep\-re\-sen\-ta\-tion in which
blue components $H_0, \dots, H_{l-1}$ are periodic of period $l$.
There is a vertex $\dot\xi_0$ in $H_0$ and a word $u_0$
so that the vertex $\dot\zeta_0$ on $H_0$ obtained by pulling back $l$ red steps
from $\dot\xi_0$ via a word $v_0$ satisfies either
\begin{align}
 \pi(f_{v_0}) \dot\zeta_0 = \dot\xi_0 &\qand
 \pi(e_{u_0}) \dot\xi_0 = \dot\zeta_0 \tag{3bi}\\
\intertext{or}
 \pi(f_{v_0}) \dot\zeta_0 &= \dot\xi_0 = \pi(e_{u_0}) \dot \zeta_0  .\tag{3bii}
\end{align}
This same type of relationship persists for each vertex obtained from
pulling back along every path leading into $\dot\xi_0$.
\end{lem}

\begin{pf}
Since $H_0$ is an infinite tail representation of $\A_m$,
there is a vertex $\dot\zeta_0$ in $H_0$ and words $u_1$ and $u_2$ so that
\[ \pi(e_{u_1}) \dot\zeta_0 = \dot\xi \qand  \pi(e_{u_2}) \dot\zeta_0 = \dot\zeta ; \]
and thus $\pi(f_v e_{u_2}) \dot\zeta_0 = \dot\xi$.
Use the commutation relations to write $f_v e_{u_2} = e_{u'_2} f_{v_0}$.
Set $\dot\xi_0 = \pi(f_{v_0}) \dot \zeta_0$.
Then
\[
 \pi(e_{u_1})\dot  \zeta_0 = \dot \xi = \pi(f_v e_{u_2}) \dot \zeta_0
 =  \pi(e_{u'_2}) \pi(f_{v_0}) \dot \zeta_0 = \pi(e_{u'_2}) \dot\xi_0.
\]
It follows that the vertices $\dot \zeta_0 $ and $\dot\xi_0$ are obtained from
$\dot\xi$ by pulling back along the blue edges $|u_1|$ and $|u'_2|=|u_2|$
steps respectively.
So they are comparable, and the relationship depends on whether
$|u_2|$ is less than, equal to or greater than $|u_1|$.

If $|u'_2| < |u_1|$, then uniqueness of the pull back along blue edges
means that $u_1 = u'_2 u_0$.  So $\pi(u_0) \dot\zeta_0 = \dot\xi_0$.
Similarly, if $|u'_2| > |u_1|$, then $u'_2 = u_1 u_0$ and
$\pi(e_{u_0}) \dot \zeta_0 = \dot \xi_0$.
In the case $|u'_2|=|u_1|$, we have $\dot\xi_0 = \dot\zeta_0$ and so
$\dot\xi_0 = \pi(f_{v_0}) \dot \xi_0$.
This is the tail by ring case, which has been excluded.

In the first case, suppose that $\dot\eta$ is any vertex in these
$l$ components which is obtained by pulling back from $\xi_0$.
That is, there is a word $w\in \Fth$ so that $\pi(w) \dot\eta = \dot\xi_0$.
There is a basis vector $\dot\zeta$ in the same component and
a word $v'$ with $|v'|=l$ so that $\pi(f_{v'}) \dot\zeta = \dot\eta$.
 From the commutation relations, there are words $v''$ of length $l$
and $w'$ of degree $d(w)$ so that $w f_{v'} = f_{v''} w'$.
Therefore,
\[ \pi(f_{v''}) \pi(w') \dot\zeta = \dot\xi_0 = \pi(f_{v_0}) \dot\zeta_0 .\]
Uniqueness implies that $v''=v_0$ and $\pi(w') \dot\zeta = \dot\zeta_0$.

Now factor $e_{u_0}w' = w'' e_{u'}$ with $|u'|=|u_0|$ and $d(w'')=d(w')$.  Then
\[  \dot\xi_0 = \pi(e_{u_0}w') \dot\zeta = \pi(w'') \pi(e_{u'})  \dot\zeta .\]
Since $w''$ is a word with the same degree, say $(s,t)$, as
$w'$ and $w$, we deduce that $\pi(e_{u'})  \dot\zeta$ is the unique
vertex obtained by pulling back from $\dot\xi_0$ by $s$ blue
and $t$ red edges, namely $\dot\eta$.  That is,
$\pi(e_{u'})  \dot\zeta = \dot\eta$; and $\dot\zeta$ lies above $\dot\eta$
in its blue component.

The other case is handled in a similar manner.
\bx\end{pf}

We now proceed as in case 3a.
Start with the basis vertex $\dot\xi_{0,0}$ in $H_0$
provided by Lemma~\ref{L:comparable}.
Pull back on both blue and red edges to determine integers $i_{s,t}$
and $j_{s,t}$ for $s,t \le 0$ and vertices $\dot\xi_{s,t}$ so that
\begin{align*}
 \pi( e_{i_{s,t}}) \dot\xi_{s-1,t} &= \dot\xi_{s,t} \qand \\
 \pi( f_{j_{s,t}}) \dot\xi_{s,t-1} &= \dot\xi_{s,t}  \qfor s, t \le 0 .
\end{align*}
The difference in case 3b is that there is periodicity.

In case 3bi, there are unique words $u_0$ of length $k$ and
$v_0$ of length $l$ so that $\pi(f_{v_0} e_{u_0}) \dot\xi_0 = \dot\xi_0$.
By Lemma~\ref{L:comparable}, for each vertex $\xi_{s,t}$, pulling
back $k$ blue edges and $l$ red edges will return to the same vertex.
That is, $\dot\xi_{s-k,t-l}=\dot\xi_{s,t}$ for all $s,t \le 0$.
So there is $\bZ(k,l)$ periodicity.
This allows us to extend the definitions to all of $\bZ^2$ using
the periodicity; and to then collapse this to a representation on
$\bZ^2/\bZ(k,l)$.
Write $\xi_{[s,t]}$ or $\xi_g$ for the vector associated to
an element $g=[s,t]$ coming from the equivalence class of $(s,t)$.
It is now easy to see that the vectors $\xi_g$ are distinct
because one can always choose the representative $(s,t)$
with $0 \le t < l$, determining the component $H_t$, and
within this component, $s$ determines the position on the
infinite tail.  Indeed, the $(k,l)$-periodicity allows us to select
a distinguished spine because pulling back $l$ steps along
the red edges moves us forward $k$ (specific) blue edges.
We shall see soon that in this case, there is always additional
symmetry.

In case 3bii, there are unique words $u_0$ of length $k$ and
$v_0$ of length $l$ so that
$\pi(f_{v_0}) \dot\zeta_0 = \pi( e_{u_0}) \dot\zeta_0 = \dot\xi_{0,0}$.
Hence $\dot\zeta_0 = \dot\xi_{0,-l}=\dot\xi_{-k,0}$.
By Lemma~\ref{L:comparable}, for each vertex $\xi_{s,t}$, pulling
back $k$ blue edges or pulling $l$ red edges will result in the same vertex.
That is, $\dot\xi_{s-k,t}=\dot\xi_{s,t-l}$ for all $s,t \le 0$.
So there is $\bZ(k,-l)$ periodicity.
In this case, there is no canonical way to carry forward.
However, as in the previous case, we obtain a
parameterization of the basis vectors as a semi-infinite
subset of $\bZ^2/\bZ(k,-l)$; and we will write
$\xi_{[s,t]}$ or $\xi_g$ for the vector associated to
an element $g=[s,t]$ coming from the equivalence class of $(s,t)$
when there is a representative with $s,t \le 0$.

In both cases, the subspace $\M = \spn\{\xi_{[s,t]} : s,t \le 0 \}$
is a coinvariant cyclic subspace.  So it is sufficient to determine
shift tail equivalence.  As usual, we define the symmetry subgroup
$H_\pi$ from the shift tail symmetry of the data
$\Sigma(\pi,\xi_{0,0}) = \{ (i_{[s,t]}, j_{[s,t]}) : s,t \le 0 \}$,
but consider it as a subgroup of $\bZ^2/\bZ(k,\pm l)$ by modding out
by the known symmetry.  Let $\Sigma(\pi)$ denote the shift tail
equivalence class of $\Sigma(\pi,\xi_{0,0})$.
The following result follows in an identical manner to the 3a case,
so it will be stated without proof.
The converse to the second statement does not follow
because there are issues with the scalars that will be dealt with soon.

\begin{thm}
If $\pi$ is a type 3b atomic $*$-rep\-re\-sen\-ta\-tion,
then $\Sigma(\pi)$ is independent of the choice of initial vector.
If two type 3b $*$-rep\-re\-sen\-ta\-tions $\pi_1$ and $\pi_2$ are unitarily equivalent,
then $\Sigma(\pi_1) = \Sigma(\pi_2)$.
\end{thm}

We now consider these two cases in more detail.

\medbreak\textbf{Case 3bi.}
Here we have symmetry $\bZ(k,l)$ with $k,l > 0$.
In this case, there are words $u_0 = i_{[-1,0]}\dots i_{[-k,0]}$
and $v_0 = j_{[0,l-1]}\dots j_{[0,0]}$ so that
\[ \pi(e_{u_0} f_{v_0})\dot\xi_{[0,0]} = \dot\xi_{[0,0]} .\]
The basis $\dot\xi_{[s,t]}$ satisfies $\dot\xi_{[ks,0]}=\dot\xi_{[0,-ls]}$.
Since we can pull back along either blue or red edges,
we find that these vectors are defined for all $s \in \bZ$.
Indeed, all $\dot\xi_{[s,t]}$ are defined for $(s,t) \in \bZ^2$.

This leads to the important observation about Case 3bi that this
representation is determined by $u_0$ and $v_0$.

\begin{lem}\label{3bi_period}
An atomic $*$-rep\-re\-sen\-ta\-tion of type {\em 3bi} is determined by
words $u_0$ and $v_0$ and a constant $\beta \in \bT$
satisfying $\pi(e_{u_0} f_{v_0})\xi_{0,0} = \beta \xi_{0,0}$.
The symmetry subgroup $H_\pi$ is always non-zero.
Indeed, there is full symmetry (without tail equivalence),
namely there is an integer $p>0$ so that
\[
 i_{[s-p,t]} = i_{[s,t]} \AnD j_{[s-p,t]} = j_{[s,t]}
 \qforal s, t \in \bZ .
\]
\end{lem}

\begin{pf}
We can find words $u_1$, $v_1$, $u_{-1}$ and $v_{-1}$ to factor
\[
 e_{u_0} f_{v_0} = f_{v_1} e_{u_1}  \qand
 f_{v_0} e_{u_0} = e_{u_{-1}} f_{v_{-1}} .
\]
Continuing recursively, we obtain words $u_r$ in $\{1,\dots,m\}^k$
and $v_r$ in $\{1,\dots,n\}^l$ for $r \in \bZ$ so that
\[ e_{u_r} f_{v_r} = f_{v_{r+1}} e_{u_{r+1}} \qforal r \in \bZ .\]
This determines the doubly infinite paths
\[
 \tau_e = \dots u_1 u_0 u_{-1} \dots \qand
 \tau_f = \dots v_1 v_0 v_{-1} \dots .
\]
The $e$'s and $f$'s move in opposite directions,
and the two spines intersect every $k$ steps forward along the blue path
for every $l$ steps backward along the red path.
The commutation relations allow us to compute
the $l$ infinite blue paths and $k$ infinite red paths,
completing the picture of this coinvariant subspace.

The commutation of words $e_u$ of length $k$ with words
$f_v$ of length $l$ is given by a permutation $\theta'$
in $S_{m^k \times n^l}$ determined by $\theta$ so that
$e_uf_v = f_{v'} e_{u'}$, where $\theta'(u,v)=(u',v')$.
The pairs $(u_r,v_r)$ therefore satisfy
$\theta'(u_r,v_r) = (u_{r+1},v_{r+1})$.
It follows that the pairs $(u_r,v_r)$ move repeatedly through a
cycle of the permutation $\theta'$.
Consequently, the sequence $(u_r,v_r)$ is periodic
of length $p$, where $p$ is the length of the cycle.
It follows that
\[ i_{[s,t]} = i_{[s+pk,t]} \qand j_{[s,t]} = j_{[s,t-pl]} = j_{[s+pk,t]} .\]
Therefore $H_\pi$ contains $[pk,0]=[0,-pl]$; and thus
is a non-trivial symmetry subgroup.
Moreover the symmetry is global, not just for $s,t \le T$.

The scalars are determined by Theorem~\ref{T:scalars}.
The character $\psi$ of $\bZ(k,l)$ is given by $\psi(k,l) = \beta$
where $\pi(e_{u_0} f_{v_0})\xi_{0,0} = \beta \xi_{0,0}$.
One can extend this to a character on $\bZ^2$; and we will
do this by setting $\phi(1,0) = 1$ and $\phi(0,1)=\beta_0$
where $\beta_0$ is an $l$th root of $\beta$.
\bx\end{pf}

As an immediate consequence of the previous two results,
we obtain:

\begin{cor}
Two type 3bi $*$-rep\-re\-sen\-ta\-tions $\pi_1$ and $\pi_2$ are unitarily equivalent
if and only if $\Sigma(\pi_1)=\Sigma(\pi_2)$ and $\beta_1 = \beta_2$,
where $\beta_1$ and $\beta_2$ are the scalars of Lemma~$\ref{3bi_period}$.
\end{cor}

Conversely, we obtain

\begin{thm}\label{T:3bi}
Given words $u_0$ in $\{1,\dots,m\}^k$
and $v_0$ in $\{1,\dots,n\}^l$ and scalar $\beta \in \bT$,
there is an atomic $*$-rep\-re\-sen\-ta\-tion of $\Fth$ of type 3bi
determined by this data.
\end{thm}

\begin{pf}
The proof of Lemma~\ref{3bi_period} explains how this is done.
The pair $(u_0,v_0)$ lies in a cycle of $\theta'$
\[ \big( (u_0,v_0), (u_1,v_1), \dots, (u_{p-1},v_{p-1}) \big) .\]
So $e_{u_i}f_{v_i} = f_{v_{i+1}} e_{u_{i+1}}$ for $i \in \C_p$.
As in the construction of Examples~\ref{E:cycle reps}
and \ref{cycle}, we obtain that $e_{u_{p-1} \dots u_1u_0}$
commutes with $f_{v_0v_1 \dots v_{p-1}}$.
Indeed, one has the useful relations
\[
 e_{u_i} f_{v_i v_{i+1} \dots v_{p-1} v_0 \dots v_{i-1}} =
 f_{v_{i+1} \dots v_{p-1} v_0 \dots v_i} e_{u_i} .
\]
This is the algebraic form of the $(k,l)$-periodicity.

By Lemma~\ref{L:ringbyring}, there is a finitely correlated
defect free representation on $\C_{pk} \times \C_{pl}$
determined by this commuting pair.
The relations above ensure $(k,l)$ periodicity;
so there is a corresponding construction on $\C_{pk} \times \C_{pl}/\ip{(k,l)}$.

The idea is to `unfold' this to obtain a 3bi $*$-rep\-re\-sen\-ta\-tion.
It is probably easier to envisage unfolding the representation
on $\C_{pk} \times \C_{pl}$ to a representation on $\bZ^2$
with $pk\bZ \times pl\bZ$ symmetry, and then observing
that the $(k,l)$-periodicity of the original picture becomes
$\bZ(k,l)$ periodicity of the type 3a representation.
So one now goes to the quotient $\bZ^2/\bZ(k,l)$ to
obtain the desired representation of type 3bi.

One can deal with scalars as before.
\bx\end{pf}

Thus we obtain:

\begin{cor}
Every $*$-rep\-re\-sen\-ta\-tion of type 3bi comes from a group construction
for a group of the form $\bZ^2/\bZ(k,l)$ with $kl>0$.
It always has a non-trivial symmetry group, and so is never
irreducible.  It decomposes as a direct integral of irreducible
$*$-rep\-re\-sen\-ta\-tions of \mbox{type 1.}
\end{cor}

\bigbreak
\medbreak\textbf{Case 3bii.}
This is the trickiest case.
Here we have $\bZ(k,-l)$ symmetry, where $k,l > 0$.
That is, our basis is $\xi_{[s,t]}$ for $s,t \le 0$ where the
equivalence class consists of cosets of $\bZ(k,-l)$.
We have a word $u_0$ of length $k$ and $v_0$ of length $l$ so that
$\pi(e_{u_0}) \dot\xi_{[-k,0]} = \pi(f_{v_0}) \dot\xi_{[-k,0]} = \dot\xi_{[0,0]}$.
Thus there is a scalar $\beta\in\bT$ so that
\[ \pi(f_{v_0}) \xi_{[-k,0]} = \beta \pi(e_{u_0}) \xi_{[-k,0]} .\]

Pulling back from $\dot\xi_{[0,0]}$ along the blue edges yields the infinite tail
\[ \tau_e = i_{[0,0]}i_{[-1,0]}i_{[-2,0]}\dots = u_0u_1u_2 \dots \]
where $u_d=i_{[dk,0]}i_{[dk-1,0]}\dots i_{[(d-1)k+1,0]}$ for $d\le0$ are
the consecutive words of length $k$.
Similarly, pulling back from $\dot\xi_{[0,0]}$ along the
red edges yields the infinite tail
\[  \tau_f = j_{[0,0]}j_{[0,-1]}j_{[0,-2]}\dots =v_0v_{-1}v_{-2}\dots \]
where $v_d=i_{[0,dl]}i_{[0,dl-1]}\dots i_{[0,(d-1)l+1]}$ for $d\le0$ are
the consecutive words of length $l$.

Let $\beta_0$ be chosen so that $\beta_0^l = \beta$.
Then as before, we may normalize the basis vectors so that
\[
 \pi(e_{i_{[s,t]}}) \xi_{[s-1,t]} = \xi_{[s,t]} \qand
 \pi(f_{j_{[s,t]}}) \xi_{[s,t-1]} = \beta_0 \xi_{[s,t]}
\]
for $sl+tk \le 0$.

Since $[dk,0]=[0,dl]$, we have
\[
 \pi(f_{v_{d+1}}e_{u_d})\dot\xi_{[dk,0]} = \dot\xi_{[(d+2)k,0]}
 = \pi(e_{u_{d+1}}f_{v_d})\dot\xi_{[dk,0]};
\]
and therefore we obtain the commutation relations
\[ f_{v_{d+1}}e_{u_d} = e_{u_{d+1}}f_{v_d} \qfor d < 0 .\]
This is a rather strong compatibility condition, and suggests why
not all pairings are possible.

For each $d \le 0$, there are at most $m^k n^l$ possible pairs $(u_d,v_d)$.
By the Pigeonhole Principle, one of these pairs
is repeated infinitely often.
Without loss of generality, we may suppose that this
sequence begins at 0; so that there are integers
$0=d_0 > d_1 > d_2 > \dots$ so that
\[
 u_{d_r} = u_0 \qand v_{d_r} = v_0 \qfor r \ge 1 .
\]
A computation now shows that $e_{u_0u_1\dots u_{d_r+1}}$ and
$f_{v_0v_{-1} \dots v_{d_r+1}}$ commute.
Indeed, one readily computes that
\[
 e_{u_i\dots u_{j-1}} f_{v_j} = f_{v_i} e_{u_{i+1}\dots u_j}
 \qand
 e_{u_i} f_{v_{i+1}\dots v_j} = f_{v_i\dots v_{j-1}} e_{u_j} .
\]
This encodes the $\bZ(k,-l)$ symmetry.
Repeated application of this yields
\begin{align*}
 e_{u_0u_1\dots u_{d_r+1}} f_{v_0v_{-1} \dots v_{d_r+1}} &=
 e_{u_0u_1\dots u_{d_r+2}}e_{u_{d_r+1}} f_{v_{d_r}} f_{v_{-1}v_{-2} \dots v_{d_r+1}} \\&=
 e_{u_0u_1\dots u_{d_r+2}}f_{v_{d_r+1}} e_{u_{d_r}} f_{v_{-1}v_{-2} \dots v_{d_r+1}} \\&=
 f_{v_0} e_{u_1\dots u_{d_r+1}} e_{u_0} f_{v_{-1}v_{-2} \dots v_{d_r+1}} \\&=
 f_{v_0} e_{u_1\dots u_{d_r+1}} f_{v_0\dots v_{d_r+2}} e_{u_{d_r+1}} \\&=\vdots\\&=
 f_{v_0v_{-1} \dots v_{d_r+1}} e_{u_0u_1\dots u_{d_r+1}} .
\end{align*}

Now we can construct a sequence of ring by ring representations
on the finite groups $G_r = \C_{|d_r|k}\times \C_{|d_r|l}/\ip{(k,-l)}$ by
building a ring by ring representation on $G_r$ using these
commuting words and the fact that they have the $\bZ(k,-l)$ symmetry.
That is, consider a finitely correlated representation $\rho_r$ on a
basis $\zeta_g$ for $g = [s,t]$ in $G_r$ for $d_rk \le s < 0$ and $0 \le t < l$ by:
\begin{align*}
 \rho_r(e_{i}) \zeta_{[s-1,t]} &= \delta_{i\,i_{s,t}}\, \zeta_{[s,t]} \qand
 \rho_r(f_{j}) \zeta_{[s,t-1]} = \delta_{j\,j_{s,t}}\, \beta_0 \zeta_{[s,t]} .
\end{align*}

These observations also provide a way (indeed many ways)
to extend the definition of the coinvariant subspace to
the full group $\bZ^2/\bZ(k,l)$ by making it $(d_1k, -d_1l)$
periodic moving forward.  That is, one may define
\[
 i_{[s,t]} = i_{[s',t']} \AnD j_{[s,t]} = i_{[s',t']}
 \qfor sk+tl > 0,
\]
where $(s',t')$ is chosen in $(d_1k,0]\times(d_1l,0]$
so that $s' \equiv s \pmod{d_1k}$ and $t'\equiv t \pmod{d_1l}$.

\begin{lem}\label{L:3bii}
An atomic $*$-rep\-re\-sen\-ta\-tion of type {\em 3bii} is determined
by a scalar $\beta \in \bT$ and
two infinite tails $\tau_e = u_0u_1u_2 \dots$ and
$\tau_f = v_0v_{-1}v_{-2}\dots$, where $|u_d|=k$ and $|v_d|=l$,
satisfying
\[ f_{v_{d+1}}e_{u_d} = e_{u_{d+1}}f_{v_d} \qfor d < 0 .\]
\end{lem}

\begin{pf}
The discussion above shows that $\beta$ and two infinite words with the
desired properties are associated to each $*$-rep\-re\-sen\-ta\-tion of type 3bii.
Conversely, suppose that such data is given.
Then by Section~\ref{S:atomicDC}, there are finitely correlated
representations $\rho_r$ as defined above for $r\ge1$.
The data $\{ (i_{s,t}, j_{s,t}) : s,t \le 0 \}$ has $\bZ(k,-l)$ symmetry,
and may be extended in a compatible way to all of $\bZ^2/\bZ(k,-l)$.
Therefore we may construct a representation on $\bZ^2/\bZ(k,-l)$
using the group construction of Section~\ref{S:group}.
\bx\end{pf}

Again we consider $\Sigma(\pi,\xi)$ to be defined on the group $\bZ^2/\bZ(k,-l)$;
and define $\Sigma(\pi)$ to be its equivalence class modulo shift tail equivalence.
The scalars are determined by Theorem~\ref{T:scalars}
by the character $\psi$ on $\bZ(k,-l)$ given by
$\pi(e_{u_0}) \pi(f_{v_0})^* \xi_{0,0} = \psi(k,-l)\xi_{0,0}$.
Combining the previous results yields:

\begin{cor}
Two type 3bii $*$-rep\-re\-sen\-ta\-tions $\pi_1$ and $\pi_2$ are unitarily equivalent
if and only if $\Sigma(\pi_1)=\Sigma(\pi_2)$ and $\psi_1 = \psi_2$, where
$\psi_i$ is the character on $\bZ(k,-l)$ determined by $\pi_i$.
\end{cor}

\medbreak\textbf{Decomposition.}
It now follows that a type 3 $*$-rep\-re\-sen\-ta\-tion has the form of one of the
group constructions.  In particular, by  Lemma~\ref{L:irreducible}, a
$*$-rep\-re\-sen\-ta\-tion $\pi$ is irreducible if and only if the symmetry
group $H_\pi$ is trivial.  In general, we obtain a direct integral
decomposition into irreducibles.
The case 3bi is never irreducible, so these representations
decompose as a direct integral of finitely correlated $*$-rep\-re\-sen\-ta\-tions.

\begin{thm}\label{Decomp3}
A (connected atomic) tail by tail $*$-rep\-re\-sen\-ta\-tion $\pi$ with symmetry group
$H_\pi \le G_\pi$ decomposes as a direct integral of irreducible
atomic $*$-rep\-re\-sen\-ta\-tions dilating a family of representations on $\ltwo(G_\pi/H_\pi)$.
\end{thm}

Consider the possibilities for case 3a, the inductive type.
The subgroup $H_\pi = \ip{(k,l)}$ for $kl>0$ cannot occur, because
this would require an irreducible $*$-rep\-re\-sen\-ta\-tion of type 3bi.

\begin{cor}\label{3asymm}
If $\pi$ is a tail by tail $*$-rep\-re\-sen\-ta\-tion of inductive type, then the symmetry
group is one of the following:
\begin{enumerate}
\item $H_\pi =\{(0,0)\}$ when $\pi$ is irreducible.
\item $H_\pi = \ip{(k,l)}$ where $kl<0$, and $\pi$ is a direct integral of irreducible
 {\em 3bii} $*$-rep\-re\-sen\-ta\-tions.
\item $H_\pi = \ip{(k,l)}$ where $kl=0$, and $\pi$ is a direct integral of irreducible
 type {\em 2} $*$-rep\-re\-sen\-ta\-tions.
\item $H_\pi$ has rank $2$, and $\pi$ is a direct integral of irreducible
 finitely correlated $*$-rep\-re\-sen\-ta\-tions.
\end{enumerate}
\end{cor}

It is natural to ask whether $\Fth$ has irreducible $*$-rep\-re\-sen\-ta\-tions of inductive type.
This is equivalent to the existence of an infinite tail without any periodicity.
This is exactly the aperiodicity condition introduced by Kumjian and Pask \cite{KumPask}.
This property is generic, but there are periodic examples such as the flip
semigroup of Example~\ref{E:flip}.
In \cite{DYperiod}, this property is explored in detail.
The semigroup $\Fth$ is either aperiodic, or it has $\bZ(k,-l)$ periodicity
for some $kl>0$.  This leads to structural differences in $\ca(\Fth)$.

It would be interesting to determine whether all aperiodic $\Fth$ have
irreducible $*$-rep\-re\-sen\-ta\-tions of type 3bii.



\end{document}